\documentclass[final,3p]{elsarticle}

\usepackage[utf8]{inputenc}
\usepackage[dvipsnames]{xcolor}
\usepackage{amsmath,amssymb,amsthm}
\usepackage{physics}
\usepackage{hyperref}
\usepackage{subcaption}
\usepackage{bbm}
\usepackage{makecell}

\usepackage{graphicx}%
\usepackage{multirow}%
\usepackage{amsmath,amssymb,amsfonts}%
\usepackage{amsthm}%
\usepackage{mathrsfs}%
\usepackage[title]{appendix}%
\usepackage[dvipsnames]{xcolor}%
\usepackage{textcomp}%
\usepackage{manyfoot}%
\usepackage{booktabs}%
\usepackage{algorithm}%
\usepackage{algorithmicx}%
\usepackage{algpseudocode}%
\usepackage{listings}%
\usepackage{pdflscape}

\usepackage[T1]{fontenc}
\usepackage[english]{babel}
\usepackage{tikz}
\usepackage{pgfplots}
\usetikzlibrary{patterns}
\pgfplotsset{
    tick label style={
        /pgf/number format/use comma=false,
        /pgf/number format/1000 sep={}
    }
}

\usepackage{xcolor}

\hypersetup{colorlinks = true, allcolors = blue}
\usepackage[nameinlink,noabbrev,capitalize]{cleveref}

\newcommand{\Real}{\mathbb{R}}

\DeclareMathOperator*{\argmin}{\mathrm{argmin}}

\newtheorem{proposition}{Proposition}[section]
\newtheorem{corollary}{Corollary}[section]
\newtheorem{lemma}{Lemma}[section]

\theoremstyle{remark}

\crefname{hypothesis}{Hypothesis}{Hypotheses}
\Crefname{hypothesis}{Hypothesis}{Hypotheses}

\input Typocaps.fd

\definecolor{customOrange}{HTML}{EB8934}
\definecolor{customTeal}{HTML}{538A96}
\definecolor{customBlue}{HTML}{34C9EB}
\definecolor{customBrown}{HTML}{6B5B4D}

\begin{document}

\begin{frontmatter}
\title{On-the-Fly Reduced Bases Through the Lens of Augmented Krylov Methods}

\author[GIREF]{Philippe-André Luneau}
\ead{philippe-andre.luneau.1@ulaval.ca}
\address[GIREF]{Groupe Interdisciplinaire de Recherche en \'El\'ements Finis de l'Universit\'e Laval, D\'epartement de Math\'ematiques et Statistique, Universit\'e Laval, Qu\'ebec, Canada}
\cortext[cor1]{Corresponding author}

\begin{abstract}

On-the-fly reduced order models (ROMs) are widely used in applications where multiple queries of large-scale models are required in an iterative method, like PDE-constrained optimization or fluid-structure interaction. When compared to classical offline--online ROMs, they have the advantage of being agnostic to any precomputed training data, mesh or geometry, which makes them more general linear/nonlinear solver accelerators. The main goal of this work is to study the effect of using a reduced basis constructed on-the-fly to accelerate a restarted Krylov algorithm, and to formally bridge this category of data-driven methods with classical augmented Krylov methods. The guiding principles of on-the-fly reduced basis will be translated into the framework of Krylov subspace methods, for a Petrov-Galerkin projection method. Previous corrections generated by the algorithm are stored in a snapshot matrix, and a singular value decomposition (SVD) is used to determine a reduced subspace containing the most significant directions. It will be shown that this approach can be embedded in the framework of Krylov methods with augmented or deflated restart, but with a different, data-driven recycling strategy. An analytical study will be presented to determine the consequences of augmentation by ROM-based vectors. Numerical experiments on a convection-diffusion problem with various degrees of convection will show that the proposed methodology yields results comparable to classical iterative linear solvers in some regimes, especially for nearly symmetric or skew-symmetric operators. At equal sizes of the augmentation basis, the convergence was sometimes faster with the proposed framework than with deflation methods.

\end{abstract}


\end{frontmatter}

\section{Introduction}\label{sec:intro}

A popular class of algorithms to solve a linear system $Ax=b$ is Krylov subspace methods, or more generally, projection methods. They are very efficient compared to direct methods, because their complexity is usually upper bounded by the number of matrix-vector products (matvecs) they require. The cost per matvec is at most $O(n^2)$ for dense matrices, or $O(\mathrm{nnz})$ for sparse matrices, where nnz is the number of nonzero entries. This makes them especially suitable for very large and sparse systems arising from the numerical solution of partial differential equations (PDEs). The other driving cost of those methods is the orthogonalization step. For full projection methods, the most recent residual has to be orthogonalized against all previously computed residuals. When the algorithm converges in a large number of iterations, which can happen when the matrix is poorly conditioned, this cost becomes a non-negligible overhead per iteration. To alleviate this phenomenon, various strategies have been implemented since the popularization of these algorithms in the 1970s  and 1980s, for instance with the arrival of MINRES \cite{Paige_Saunders_1975} and GMRES \cite{Saad_Schultz_1986}. One possible approach is to keep a reduced window of previous vectors to orthogonalize against, like the so-called Incomplete Orthogonalization Method (IOM) does \cite{Brown_Hindmarsh_1986,Saad_2003}. The orthogonalization cost is then capped by the fixed size of the window. Another strategy is restarting, which consists of prematurely stopping the algorithm past a certain number of iterations, and then restarting it from the current partial solution, while completely resetting the history of vectors to orthogonalize against. While this reduces the work dedicated to orthogonalization per iteration, the convergence is drastically slowed down since components that were constructed in previous cycles are forgotten and must be recomputed. The idea of recycling and augmenting came about to counter this phenomenon. By keeping only a few vectors that summarize well the work done in a previous cycle and augmenting the search space with those, the essential information is carried forward and the convergence becomes closer to that of the full, untruncated projection method. There are multiple ways to select such information, which give rise to multiple restarted Krylov methods with augmentation or recycling. We refer to the work of Soodhalter, de Sturler and Kilmer \cite{Soodhalter_de_Sturler_Kilmer_2020} for a comprehensive review of subspace recycling. Eiermann, Ernst and Schneider \cite{Eiermann_Ernst_Schneider_2000} classify augmentation methods into two fundamental strategies. One is based on a geometrical argument of measuring largest principal angles between two subspaces to preserve the most important orthogonality constraints, which amounts to the GCROT method \cite{de_Sturler_1999}. They determine how to optimally truncate the search space of the previous cycle according to this metric.
The other strategy is based on the idea of deflation, which means that information about invariant subspaces is carried over from cycle to cycle, because it usually slows down the convergence. This information is usually obtained by approximating eigenvectors by Ritz vectors, and augmenting the search space basis with those. This approach led to many algorithms like GMRES-DR \cite{Morgan_2002}, GCRO-DR \cite{Parks_de_Sturler_Mackey_Johnson_Maiti_2006} and Deflated-CG \cite{Saad_Yeung_Erhel_Guyomarc’h_2000}. Deflation techniques are also used to design preconditioners \cite{Frank_Vuik_2001,Nabben_Vuik_2006,Erlangga_Nabben_2008} for full projection methods.

It is known that eigenvectors (or approximate eigenvectors) are a good basis for a deflation subspace because they reduce the effective condition number of the matrix by eliminating specific eigenvalues from the problem. However, it is fair to consider augmenting the basis with other vectors that are not necessarily approximations of eigenvectors, but may carry important information about the convergence. For instance, Baker, Jessup and Manteuffel \cite{Baker_Jessup_Manteuffel_2005} proposed the LGMRES method, which stores approximate errors between the cycles and augments the search space of GMRES with them. More recently, data-driven methods have been introduced. They usually reuse information obtained from a sequence of previously solved, similar linear systems in order to create augmentation spaces, especially with the proper orthogonal decomposition (POD).
Kerfriden, Gosselet, Adhikari and Bordas \cite{Kerfriden_Gosselet_Adhikari_Bordas_2011} used a POD-augmented preconditioned conjugate gradient (PCG) method in a nonlinear Newton-Krylov solver. 
Carlberg, Forstall and Tuminaro \cite{Carlberg_Forstall_Tuminaro_2016} proposed a goal-oriented POD, for instance by introducing weights, to improve the convergence of PCG when using a POD-augmented search space. Diaz Cortes, Vuik and Jansen \cite{Diaz_Cortes_Vuik_Jansen_2018} proposed the use of POD deflation vectors in the case when the right-hand sides of the linear systems are subject to low-dimensional perturbations. Their follow-up work \cite{Diaz_Cortés_Vuik_Jansen_2021} compares an on-the-fly in time approach to build the POD deflation vectors with the classical offline-online paradigm for POD \cite{Quarteroni_Manzoni_Negri_2016}. Li, Zikatanov and Zuo \cite{Li_Zikatanov_Zuo_2025} proposed the Reduced Krylov Basis Method (RKBM) framework, in which a reduced Krylov basis is constructed in the offline POD training phase, and then used to find solutions quicker during the online phase. More recently, machine learning approaches, like DeepONet, have been used as well to construct the deflation space \cite{Kopaničáková_Lee_Karniadakis_2025}.

This article proposes an approach using information collected on-the-fly in a restarted Krylov solver, to build a reduced basis following the principles of the POD framework: constructing a snapshot matrix of corrections in the Krylov subspace, and using the singular value decomposition (SVD) to retain a low-dimensional, dominant subspace spanned by those corrections. The search space for corrections in the following iterations is then initially limited to that reduced space. We can show that this method acts formally as a restarted Krylov method with a POD-augmented subspace that changes at every cycle. However, it differs from other methods of this type present in the literature because the information necessary for the construction of the reduced basis is not obtained from the solution of a sequence of different linear systems, but purely from the convergence history of the current solver. In fact, this makes this method more similar to classical augmented or deflation-restarted Krylov methods, just like LGMRES, GCROT or GMRES-DR, the difference lying in the way the augmentation space is selected. The proposed approach uses a data-driven selection criterion, like POD, instead of the usual recency, spectral or geometric selection criteria.

\Cref{sec:prelim} will recall useful concepts from the theory of Krylov subspace methods, introduce augmented Krylov methods, and present reduced bases in the on-the-fly paradigm. In \Cref{sec:theory}, reduced-basis-augmented minimal residual methods will be presented. The error between those augmented, restarted variants and untruncated methods will be derived explicitly. Moreover, we are able to characterize the proposed method as an interpolation between two classical methods, LGMRES($m$, 1) and $s$-ORTHOMIN(1) \cite{Chronopoulos_1991}. \Cref{sec:numerical} will provide numerical illustration of the capabilities of the proposed algorithms. The goal of this section is to demonstrate that, while they neither systematically nor significantly outperform existing augmented Krylov methods, they are competitive when compared to other state-of-the-art linear solvers which are not data-driven.

\section{Preliminaries}\label{sec:prelim}

\subsection{Projection Methods for Linear Systems}

Let the linear system $Ax=b$, with $A\in\Real^{n\times n}$ nonsingular and $x,b\in\Real^n$ have the solution $x=x^*$. Starting from an initial guess $x^0$, the core idea of projection methods is to look for an approximate solution of the form $x^* \approx x^k \in x^0 + \mathcal{W}_k$, with the additional property $r^k \perp \mathcal{V}_k$, where $\mathcal{W}_k,\mathcal{V}_k$ are $k$-dimensional subspaces and $r^k = b - Ax^k$ is the residual \cite{Saad_2003}. The solution is found in an affine search space, under the constraint that its residual must be orthogonal to some linear subspace. Projection methods can be classified according to the choice of the subspaces.

Petrov-Galerkin methods are the most general. Assuming only that $W_k$, $V_k$ are orthonormal bases of $\mathcal{W}_k$ and $\mathcal{V}_k$ respectively, the two conditions below must be satisfied:
\begin{align}
\label{eq:PGsystem}
    x^k = x^0 + W_k y^k \qq{and} V_k^\top (b - Ax^k) = 0 \iff V_k^\top A W_k y^k = V_k^\top r^0.
\end{align}
In general, for $A$ invertible, it is not guaranteed that the projected system even has a unique solution. Very generic hypotheses like $\mathrm{im}(AW_k) \cap \mathrm{im}(V_k)^\perp = \{0\}$ can be considered, but they are usually not very informative in practical cases. Within Petrov-Galerkin methods, the choice $\mathcal{V}_k = M\mathcal{W}_k$ yields most methods used in practice \cite{Gaul_Gutknecht_Liesen_Nabben_2013}.

Setting $M=I$, or $\mathcal{W}_k = \mathcal{V}_k$, yields Galerkin projection methods. The approximate solution satisfies the analog of the projected system \eqref{eq:PGsystem} with $V_k = W_k$,  
\begin{align*}
     W_k^\top A W_k y^k = W_k^\top r^0.
\end{align*}
This class of methods is especially adapted for the case $A$ symmetric and positive definite (SPD), since the projected system above will also be nonsingular and SPD. Moreover, the projected system to solve is the stationarity condition for the minimization problem
\begin{align*}
    x^k = \argmin_{x \in x^0 + \mathcal{W}_k} \| x - x^* \|_A^2,
\end{align*}
the minimization of the $A$-norm error. Another common choice is setting $M=A$, or $\mathcal{V}_k = A\mathcal{W}_k$. The Petrov-Galerkin system \eqref{eq:PGsystem} becomes
\begin{align*}
    W_k^\top A^\top A W_k y^k = W_k^\top A^\top r^0.
\end{align*}
It is easy to recognize the normal equations, which are equivalent to solving the minimization problem
\begin{align*}
    \| r^k\| = \min_{y\in \Real^k} \|r^0 - AW_ky \| = \min_{x \in x^0 + \mathcal{W}_k} \| b - Ax\| = \min_{x \in x^0 + \mathcal{W}_k} \| x - x^*\|_{A^\top A}.
\end{align*}
This is why those methods are generally referred to as \textit{minimal residual} projection methods. Notice the equivalence with minimizing the error in the $A^\top A$-norm, which makes this norm (and the $A$-norm for Galerkin methods) the natural choice for measuring convergence.

Krylov methods are special cases of both Petrov-Galerkin and Galerkin methods, where in addition of imposing the orthogonality constraint, the affine subspace is restricted to be the Krylov subspace of dimension $k$,
\begin{align*}
    \mathcal{W}_k = \mathcal{K}_k(A,r^0) = \mathrm{span}\{r^0, Ar^0, A^2 r^0, \dots, A^{k-1}r^0\}.
\end{align*}
This choice is particularly useful because the inverse of the matrix $A$ can be expressed as a polynomial in $A$, $p_m(A) = A^{-1}$, of degree at most $m\leq n$ \cite{Saad_2003}. A basis of the space of polynomials of degree at most $m$ would then be able to represent exactly the action of $A^{-1}$ by computing only powers of $A$ (or even, only matvecs). Indeed, the matrix $W_{m+1}$, basis of $\mathcal{W}_{m+1}$, has columns of the form $q_i(A)r^0$  with $\mathrm{deg}(q_i) = m$. They form a basis of the space of degree $m$ polynomials in $A$. For some $y^{m+1}\in\Real^{m+1}$,
\begin{align*}
    W_{m+1}y^{m+1} = \bigg(\sum_{i=1}^{m+1} y_i^{m+1}q_i(A) \bigg) r^0.
\end{align*}
By setting $y^{m+1}$ to be the coefficients of $p_m$ in the basis $\{q_i\}_{i=1}^{m+1}$, $W_{m+1}y^{m+1} = A^{-1} r^0 = x^* - x^0$, the exact correction needed to reach the solution from the starting point $x^0$. Krylov methods are built upon this principle, by constructing iteratively an orthonormal basis for $\mathcal{K}_k(A,r^0)$ until convergence is reached. Minimal residual methods will find the coefficients by minimizing the residual over this polynomial space at each iteration, and the Petrov-Galerkin condition ensures that the residual has no component in the $A$-image of the current Krylov space, i.e., that the next coefficient will belong to the next power of $A$. Numerically, convergence in $m$ or $n$ iterations does not need to be reached exactly, as we are interested in an approximate solution; a small residual would suffice to stop the algorithm. Moreover, floating-point arithmetic prevents reaching the exact solution and might propagate errors that destabilize convergence.

Both Petrov-Galerkin and Galerkin methods implement this paradigm. Note that in exact arithmetic, as long as the two spaces chosen are equal, the algorithms will produce identical iterates. However, floating-point arithmetic allows for different implementations of the same theoretical method to give different results, with potentially different properties. For instance, some representatives of Petrov-Galerkin, or more precisely minimal-residual Krylov methods are GMRES, MINRES and GCR. The conjugate gradient (CG) method and the full orthogonalization method (FOM) are representatives of Krylov methods from the Galerkin family. Frequently, these methods differ by optimizing some part of the orthogonalization process to specific matrix types, like using the Lanczos process instead of the Arnoldi process when $A$ is symmetric. This is the main distinction between FOM and CG \cite{Saad_2003}, meaning they are otherwise equivalent in exact arithmetic. The families and methods named here are not exhaustive, and other choices of spaces can yield useful projection methods, like $\mathcal{W}_k=\mathcal{K}_k(A, r^0)$ and $\mathcal{V}_k=\mathcal{K}_k(A^\top, r^0)$ for the BiCG method.

\subsection{Augmented Krylov Methods}
The methods are based on a pattern of outer iterations or cycles (index $i\leq N_O$) and inner iterations (index $k\leq N_I$). A first group of augmented Krylov methods are projection methods where the search space has the form $\mathcal{W}_{k,i} = \mathcal{V}_i \oplus \mathcal{K}_k(A, r^i)$, meaning that it is the direct sum of a Krylov space and some other subspace with suitable properties. One of the main goals of these methods is to accelerate restarted Krylov methods. For instance, in GMRES($N_I$), which is restarted every $N_I$ iterations, the loss of orthogonality with directions obtained in previous cycles slows down its convergence considerably. Augmenting the search space with previous corrections acts as some form of memory mechanism which prevents the algorithm from generating contributions in subspaces that were already considered in the past. This approach has been implemented in LGMRES \cite{Baker_Jessup_Manteuffel_2005}, where a few (usually, fewer than 5) corrections from past cycles are added in $\mathcal{V}_i$, which effectively reduces the oscillatory behaviour appearing in the convergence of GMRES($N_I$). In that case, if cycles are seen as independent restarts of the method, LGMRES is composed of a sequence of independent GMRES calls with an augmentation subspace $\mathcal{V}_i$, which changes only between cycles. Other methods like GMRES-E \cite{Morgan_1995} use this form, but with approximate eigenvectors instead. We will refer to such algorithms as augmented FGMRES methods, because Saad \cite{Saad_1997} introduced a generic modified Arnoldi algorithm following this form. \Cref{alg:augmentedpolynomialprecond} presents this framework, more precisely the $i$-th cycle, in which a Krylov subspace $K_{N_I}(A,r^i)$ is constructed while maintaining orthogonality to $m$ additional vectors.

\begin{algorithm}
\begin{algorithmic}[1]
\State $r^i = b- Ax^i$, $z^0 = r^i/\| r^i\|$, $V_i = [v^{1,i},\dots,v^{m,i}]$
\State $[z^0,\dots,z^{N_I+m}],H_i \leftarrow \mathrm{Arnoldi}(\{z^0, Az^0,A^2z^0,\dots,A^{N_I}z^0,Av^{1,i},\dots,Av^{m,i}\})$
\State $W_i = [z^0,\dots,z^{N_I-1},v^{1,i},\dots,v^{m,i}]$
\State $Z_i = [z^0,\dots,z^{N_I+m}]$, s.t. $AW_i = Z_iH_i$
\State $\zeta^i = \argmin_{\zeta \in\Real^{N_I+m}}\| Z_i^\top r^i - H_i\zeta \|$
\State $x^{i+1} = x^i + W_i\zeta^i$
\end{algorithmic}
    
    \caption{Augmented FGMRES-Like Algorithm ($i$-th cycle)}
    \label{alg:augmentedpolynomialprecond}
\end{algorithm}

Another approach suggests that the operator of the Krylov subspace changes. Subspaces of the form $\mathcal{W}_{k,i} = \mathcal{V}_i \oplus \mathcal{K}_k(P_i A, r^i)$ are considered. This means that between each cycle, the augmented basis is updated (e.g., using deflation vectors), but also the projector $P_i$, to improve the convergence of the Krylov method within the cycle by potentially removing problematic eigenvectors. Such methods include GCROT, GCRO-DR and Deflated-CG, the last two using deflation \cite{Frank_Vuik_2001} to augment the space and to accelerate convergence. Deflation has also been used with CG \cite{Saad_Yeung_Erhel_Guyomarc’h_2000}, but the augmentation basis has to come from the solution of previous linear systems, since CG need not be restarted in general. Hicken and Zingg \cite{Hicken_Zingg_2010} proposed a variant of GCROT, GCROT($m$,$k$), which simply keeps the $k$ most recent corrections as the augmentation vectors, like LGMRES would, but also projects the operator, something absent from LGMRES. \Cref{alg:augmentedgcromethod} presents the $i$-th cycle of a GCRO-like algorithm, which uses $m$ augmentation vectors. Note that the projected operator $P_iA$ makes it so that all Arnoldi vectors are orthogonal to $\mathrm{im}(A V_i)$.

For an in-depth analysis of convergence for augmented methods, we recommend the thorough analysis of Saad \cite{Saad_1997}. We also refer to the general framework introduced by Gaul, Gutknecht, Liesen and Nabben \cite{Gaul_Gutknecht_Liesen_Nabben_2013}, which allows one to describe all augmented, deflation-like methods in a single formalism \cite{Soodhalter_de_Sturler_Kilmer_2020}. 

\begin{algorithm}
\begin{algorithmic}[1]
\State $r^i = b- Ax^i$, $V_i = [v^{1,i},\dots,v^{m,i}]$
\State $Q_iR_i = AV_i$, $C_i = Q_i$, $U_i = V_iR_i^{-1}$
\State $P_i = I - C_iC_i^\top$
\State $\tilde{r}^i = P_i r^i$, $z^0 = \tilde{r}^i/\| \tilde{r}^i\|$
\State $[z^0,\dots,z^{N_I}], H_i \leftarrow \mathrm{Arnoldi}(\{z^0, P_iAz^0,(P_iA)^2z^0,\dots,(P_iA)^{N_I}z^0\})$
\State $W_i = [U_i,z^0,\dots,z^{N_I-1}]$
\State $Z_i = [C_i,z^0,\dots,z^{N_I}]$, s.t. $AW_i = Z_iH_i$
\State $\zeta^i = \argmin_{\zeta\in\Real^{N_I+m}}\| Z_i^\top r^i - H_i\zeta\|$
\State $x^{i+1} = x^i + W_i\zeta^i$
\end{algorithmic}
    
    \caption{Augmented GCRO-Like Algorithm ($i$-th cycle)}
    \label{alg:augmentedgcromethod}
\end{algorithm}

\subsection{On-the-Fly Reduced Basis}

The on-the-fly framework relies on the collection of so-called \textit{snapshots} of a vector quantity $x^k\in\Real^n$ that is computed at each iteration of an iterative process, which are stored in the columns of a snapshot matrix $S_{N_b} \in {\Real^{n\times N_b}}$, assumed to be full-rank. We note $N_b$ the number of snapshots. Often, the snapshot matrix can also be centered, denoted $\overline{S}_{N_b}$, by subtracting the mean $\overline{x}^{N_b}$ of all snapshots from each column. Notice that this reduces the rank of the matrix to $N_b-1$. We consider a concrete application case. In a nonlinear solver, if a linear system $A_k x^k = f^k$ is solved at iterations $k=1,\dots,N_b$, then at iteration $N_b+1$, the solution space of the system is restricted to either the full span of the snapshots, $\mathrm{im}(\overline{S}_{N_b})$, or to a subspace $\mathcal{V}_{N_b}\subseteq\mathrm{im}(\overline{S}_{N_b})$. In the second case, the search space is further reduced from dimension $N_b$ to dimension $m\leq N_b \ll n$. Suppose that $\mathcal{V}_{N_b}$ is spanned by an orthonormal basis, composed of the columns of a matrix $V_{N_b}\in \Real^{n\times m}$. The solution can then be written as $x_\mathrm{RB}^{N_b+1} = \overline{x}^{N_b} + V_{N_b} z^{N_b + 1}$. There are multiple ways to determine the coefficients $z^{N_b+1}$ of the reduced solution, for instance residual minimization or Galerkin projection. In the case where there is no truncation, i.e., $\mathcal{V}_{N_b} = \mathrm{im}(\overline{S}_{N_b})$, then this method acts almost like a vector extrapolation algorithm, like RRE or MPE \cite{Brezinski_Redivo-Zaglia_Saad_2018}. The main difference is that vector extrapolation would require solving the large-scale system at every iteration to construct the next extrapolated iterate, while the goal of the reduced basis is to bypass the solution of the system as often as possible. Thus, the snapshot matrix of the vector extrapolation method would be refreshed at each iteration with a new \textit{exact} solution (often in a sliding window/FIFO pattern). On the other hand, the snapshot matrix of the on-the-fly reduced basis method would only be refreshed and recompressed every $N_I$ iterations, which can be hardcoded, but ideally should be determined dynamically by using an error estimator for the reduced solution. This distinguishes them respectively as a convergence acceleration method and as a dimensionality reduction technique.

When the truncated case is considered, i.e., $\mathcal{V}_{N_b} \subsetneq \mathrm{im}(\overline{S}_{N_b})$ and $\mathrm{dim}(\mathcal{V}_{N_b}) = m < N_b$, the SVD is often used to determine the best low-rank approximation of $\overline{S}_{N_b}$ and determine its range. Indeed, for $\overline{S}_{N_b} = U_{N_b}\Sigma_{N_b} T_{N_b}^\top$, the basis matrix $V_{N_b}$ can be selected as $V_{N_b} = [U_{N_b}]_m$  ($[\,\cdot\,]_j$ will be used to denote the restriction of a matrix to its first $j$ columns), since they are associated to the $m$ largest singular values of $\overline{S}_{N_b}$. $m$ can be chosen dynamically by ensuring a certain ratio of energy is conserved, like taking the smallest $m$ such that
\begin{align}
\label{eq:energy}
    \frac{\sum_{i=1}^m \sigma_i^2}{\sum_{i=1}^{N_b-1} \sigma_i^2} > 1 - \varepsilon_\mathrm{rb}^2. 
\end{align}
The classical way to find the reduced coefficients in the POD framework is by using a Galerkin projection, i.e., by solving the system
 \begin{align*}
     V_{N_b}^\top A_{k} V_{N_b} z^{k+1} = V_{N_b}^\top(f^{k} - A_k\overline{x}^{N_b})\quad ( N_b+1 \leq k \leq N_I),
 \end{align*}
 which ensures the residual orthogonality condition $f^k - A_k x_\mathrm{RB}^{k+1} \perp \mathrm{im}(V_{N_b})$.

 When updating the snapshot matrix, the SVD has to be recomputed to take into account the new snapshots. Since the two operations done to $\overline{S}_k$ are a mean-centering update and a column update, they can be represented as an additive rank-2 update. Thus, the full SVD of the matrix need not be recomputed, and the incremental SVD algorithm proposed by Brand \cite{Brand_2006} can be used to update the matrices $U_k$, $T_k$ and $\Sigma_k$ at a lower cost, especially for larger $m$ values. 

 On-the-fly reduced bases have been applied to accelerate a wide variety of PDE-constrained optimization algorithms, including trust-region methods \cite{Yano_Huang_Zahr_2021,Wen_Zahr_2023}, line-search methods \cite{Grundvig_Heinkenschloss_2024}, and MMA \cite{Xiao_Ma_Lu_Raghavan_Zhang_2022} (introduced by Svanberg \cite{Svanberg_1987}), as well as more generic nonlinear solvers like Newton-Krylov methods \cite{Kerfriden_Gosselet_Adhikari_Bordas_2011} or fixed-point methods \cite{Vierendeels_Lanoye_Degroote_Verdonck_2007} in fluid-structure interaction.
 
\section{Theoretical Framework}
\label{sec:theory}

\subsection{Minimal Residual Method on the Reduced Subspace}

Augmentation of Krylov methods with POD vectors is not a new concept in itself \cite{Carlberg_Forstall_Tuminaro_2016,Diaz_Cortes_Vuik_Jansen_2018,Diaz_Cortés_Vuik_Jansen_2021,Li_Zikatanov_Zuo_2025,Vallot_Vuyst_Dairay_Mougeot_2026}. However, most methods rely on parametric POD, meaning that the linear system comes from a parametrized PDE, and the solution to this PDE is computed for multiple parameter values. Then, a reduced basis is computed with the SVD, which can then be used as augmentation vectors for a Krylov solver, when solving the system at a new parameter value. However, we will leverage the POD framework in a different way, to compress information within a single linear solve. An on-the-fly reduced basis will be constructed from the convergence history. The reduced basis plays the role of an augmentation of the Krylov subspace. Indeed, the residual minimization problem
\begin{align*}
    \min_{\zeta\in\Real^m} \| r^i - AV_i \zeta\|
\end{align*}
can be seen as restricting the search space to a reduced subspace, as well as the first intermediate minimization in an augmented Krylov method.

In this section, the notation $u_\circ$ will be used to denote a zero-padded vector $u$. We propose a variant of LGMRES($N_I$, $k$) named GMRES-RB($N_b$, $N_I$, $\varepsilon_\mathrm{rb}$, $m_\mathrm{max}$). This method will run $N_b$ regular GMRES iterations, producing intermediate iterates
\begin{align*}
    x^k = x^0 + K_ky^k,\quad k=1,\dots,N_b,
\end{align*}
where $K_k\in \Real^{n\times k}$ is the Krylov subspace matrix at iteration $k$ and $y^k\in\Real^k $ is the correction in Krylov coordinates. Initially, a snapshot matrix $S_0$ will be made with the $N_b-1$ zero-padded correction vectors, $y^k_\circ \in \Real^{N_b}$, from which the final correction $y^{N_b}$ is subtracted. With an SVD $S_0=U_0\Sigma_0 T_0^\top$, the $m_0$ principal left singular vectors are extracted, by using the criterion that $100\cdot(1-\varepsilon_\mathrm{rb}^2)\%$ of the energy is captured. Because of the subtraction by $y^{N_b}$, the columns of the snapshot matrix represent the intermediate corrections as perturbations of the final correction. Heuristically, we hope this will help to detect directions that were harder to converge during a single cycle. Then, \Cref{alg:augmentedpolynomialprecond} is started with the augmentation matrix $V_0=K_{N_b}\big[[U_0]_{m_0}, y^{N_b}_\circ\big]$, which effectively contains the total correction of the GMRES cycle and a subspace spanning significant corrections. Then, the $N_I$ inner corrections associated to the inner Arnoldi iterations from $z^0$, $\zeta^{k,0}$, are stored, which is not usually the case in a typical application of \Cref{alg:augmentedpolynomialprecond}. At the end of this cycle, the POD methodology is again applied to the corrections: the last correction is subtracted from the others, and an energy-based criterion is applied to the SVD to select dominant modes. Moreover, we add a maximal reduced basis size $m_\mathrm{max}$, which can be set to $\infty$ if no threshold needs to be imposed. The next cycle of \Cref{alg:augmentedpolynomialprecond} is started with the resulting reduced basis $V_{i+1}$. A pseudocode description is available in \Cref{alg:gmresaccelere}.

\begin{algorithm}
\begin{algorithmic}[1]
\State $x_\textrm{GMRES}^0,\varepsilon_\mathrm{rb},m_\mathrm{max},N_b,\varepsilon,N_O,N_I$ given
\State Do $N_b$ GMRES iterations, produces Krylov basis matrix $K_{N_b}$ and corrections $y^i \in \Real^{i}$
\State $S_0=[y^1_\circ-y^{N_b}, \dots,y^{N_b-1}_\circ-y^{N_b}]$
\State $U_0\Sigma_0 T_0^\top = S_0$
\State Choose the smallest $m_0$ such that $[U_0]_{m_0}$ captures $100\cdot(1-\varepsilon_\mathrm{rb}^2)\%$ of the energy
\State $m_0 \leftarrow \min(m_0,m_\mathrm{max})$
\State $V_0=K_{N_b}\big[[U_0]_{m_0}, y^{N_b}\big]$
\State  $x^0_\textrm{GMRES-RB} = x^{N_b}_\textrm{GMRES}$
\State $i=0$
\While{$\|r^i\| > \varepsilon$ and $i<N_O$}
\State Apply \Cref{alg:augmentedpolynomialprecond} with augmentation vectors $V_i$ and initial data $x^i$
\For{$k=1,\dots,N_I-1$}
\State Obtain the $k$-th partial correction: $$\zeta^{k,i} = \argmin_{\zeta\in\Real^{k+m_i+1}} \|r^i - AW^{(k)}_i\zeta \|,$$where $W^{(k)}_i=[z^0,\dots,z^{k-1},V_i]$
\EndFor
\State  $ S_i = [\zeta_\circ^{1,i} -\zeta^{i},\dots,\zeta_\circ^{N_I-1,i}-\zeta^{i}]$
    \State $U_{i+1}\Sigma_{i+1} T_{i+1}^\top = S_i$
    \State Choose the smallest $m_{i+1}$ such that $[U_{i+1}]_{m_{i+1}}$ captures $100\cdot(1-\varepsilon_\mathrm{rb}^2)\%$ of the energy
    \State $m_{i+1}\leftarrow \min(m_{i+1},m_\mathrm{max})$
    \State $V_{i+1} = W_i\big[[U_{i+1}]_{m_{i+1}},\zeta^{i}\big]$
    \State $i\leftarrow i+1$
\EndWhile
\end{algorithmic}
    
    \caption{GMRES-RB($N_b,N_I,\varepsilon_\mathrm{rb},m_\mathrm{max}$)}
    \label{alg:gmresaccelere}
\end{algorithm}
 In practice, the partial residual minimization problems can be solved more efficiently by using directly restrictions of $H_i$. Also, about the padded snapshot vectors, notice that the zeros will have to be added in the middle of the vectors, to eliminate the columns from $W_i$ associated to $z^k,\dots,z^{N_I-1}$.

The present algorithm is just a wrapper around \Cref{alg:augmentedpolynomialprecond} with an SVD-based choice for the augmentation vectors. Notice that all vectors in the augmentation basis of cycle $i+1$ originate from cycle $i$; this is one of the main differences with LGMRES and GCROT($m$,$k$), where the augmentation vectors are the final corrections from the $k$ previous cycles. Also, for both of these methods, the criterion to determine which vector to keep is recency (first in, first out). Here, we propose that the final correction is always kept, but that additional vectors could be retained as well, using the energy-based criterion \eqref{eq:energy}. Moreover, GCROT($m$,$k$) is a special case of \Cref{alg:augmentedgcromethod}, which projects the operator on the orthogonal complement of the augmentation basis at each cycle.

We also introduce four other variants of the algorithm:
\begin{enumerate}
    \item GMRES-RB-O($N_b$, $N_I$, $k$), which keeps the last correction and the $k-1$ \textit{oldest} snapshots in $S_i$;
    \item GMRES-RB-L($N_b$, $N_I$, $k$), which keeps the last correction and the $k-1$ \textit{latest} snapshots in $S_i$;
    \item GMRES-RB-C($N_b$, $N_I$, $\varepsilon_\mathrm{rb}$, $m_\mathrm{max}$,  $\omega$), which compresses the \textit{cumulative} history of retained modes with the energy criterion, weighted by $\omega\in(0,1)$;
    \item GCRO-RB($N_b$, $N_I$, $\varepsilon_\mathrm{rb}$, $m_\mathrm{max}$), which uses \Cref{alg:augmentedgcromethod} as an inner loop.
\end{enumerate}

The third variant is worth exploring in more detail. In \Cref{alg:gmresaccelere}, all the recycling information comes exclusively from the previous cycle. For some harder problems, history from previous cycle might help convergence as well. This is done by realizing that the columns $N_I+1$ to $N_I+m_i+1$ of $W_i$ are $V_i = W_{i-1}[[U_i]_{m_i},\zeta^{i-1}]$, precisely the history that we want to preserve. We define
\begin{align*}
    \Gamma_i = R_{m_i} [\Sigma_i]_{m_im_i},
\end{align*}
where $[\Sigma_i]_{m_im_i}$ is the top-left, square $m_i\times m_i$ submatrix of $\Sigma_i$, and $R_{m_i}\in\Real^{(N_I+m_i+1) \times m_i}$ is a zero matrix with an $m_i\times m_i$ identity matrix in the bottom. Its columns are the coordinates of $V_i$ in the $W_i$ basis, weighted by their original energy. We can consider the following augmented snapshot matrix for the SVD,
\begin{align*}
    S_i = [ \omega\Gamma_i, \zeta_\circ^{1,i} -\zeta^{i},\dots,\zeta_\circ^{N_I-1,i}-\zeta^{i}].
\end{align*}
Notice that information that is allowed to traverse multiple cycles will decrease geometrically with rate $\omega$.

The main advantage of doing the SVD on the corrections in Krylov coordinates is that the matrix to factor is much smaller than a matrix with snapshots of $x^k$ directly.
However, since the matrix $W_i$ does not have orthonormal columns, doing the SVD on $S_i$ is not equivalent, up to a rotation, to doing it on $W_iS_i$. Experimentally, it was noticed that doing the SVD on $S_i$ directly yielded better results.

\subsection{Relations to Known Methods}

The GMRES-RB method, in exact arithmetic, can be seen as interpolating two other classical methods. In the extreme case where the energy criterion were to discard all modes ($m_\mathrm{max}=0$) and $N_b=N_I$, the method would coincide in exact arithmetic with LGMRES($N_I$, $1$); only the total correction from the previous cycle is kept. In GMRES-RB, since the augmented basis always contains the total correction of the last cycle as in LGMRES, regardless of the truncation, we know that the skip angles satisfy
\begin{align*}
    \cos\theta(r^{i+1},r^{i-1}) = \frac{\| r^{i+1}\|}{\| r^{i-1}\|},
\end{align*}
which helps to alleviate slow, alternating patterns in the corrections \cite{Baker_Jessup_Manteuffel_2005}. On the contrary, when the energy criterion does not truncate anything ($\varepsilon_\mathrm{rb}=0$ and $m_\mathrm{max}=\infty$), the method is equivalent to $N_I$-ORTHOMIN($1$), in exact arithmetic.

$s$-ORTHOMIN($k$) is an $s$-step Krylov method. These methods were first introduced by Chronopoulos \cite{Chronopoulos_1987}, with the goal of obtaining computational gains on parallel computer architectures \cite{Chronopoulos_Kim_2020}. Recent advances extended the $s$-step methods to dynamically chosen $s$ parameter, to improve convergence \cite{Shokri_Kaveh_Hajarian_Chronopoulos_2026}.  Within the $i$-th cycle, $s$-ORTHOMIN($k$) produces the iterate $x^{i+1}$ and the residual $r^{i+1}$, by residual minimization over the affine space $x^i + \mathrm{im}(P_i)$. Then, a basis of $K_s(A,r^{i+1})$ is computed, which is orthogonalized against the $k$ previous bases $P_i,\dots,P_{i-k+1}$, resulting in the new basis $P_{i+1}$. $P_0$ is set to be a basis of $K_s(A,r^0)$. It acts as a batch Krylov method, each iteration using a correction in an $s$-dimensional Krylov space. The special case $N_I$-ORTHOMIN($1$) will be detailed in \Cref{alg:orthomin(1)}.
\begin{algorithm}
\begin{algorithmic}[1]
\State $x^0$, $r^0 = b - Ax^0$, $P_0 = [r^0, Ar^0,\dots,A^{N_I-1}r^0]$
\State $i=0$
\While{$\| r^i\|>\varepsilon$ and $i<N_O$}
\State $\zeta^{i+1} = \argmin_{\zeta \in \Real^{N_I}} \| r^i - AP_i \zeta \|$
\State $r^{i+1} = r^i - AP_i\zeta^{i+1}$
\State $x^{i+1} = x^i + P_i\zeta^{i+1}$
\State $R_{i+1} = [r^{i+1}, Ar^{i+1},\dots,A^{N_I-1}r^{i+1}]$
\State $A^\top A$-orthogonalize $R_{i+1}$ against $P_i$, yields $P_{i+1}$
\EndWhile
\end{algorithmic}
    \caption{$N_I$-ORTHOMIN(1)}
    \label{alg:orthomin(1)}
\end{algorithm}

Line 4 minimizes the residual over the correction search space $\mathrm{im}(P_i)$. However, if we consider the $P_i$'s to be augmentation spaces, we would be interested in the minimization over $\mathrm{im}(P_{i-1})\oplus \mathrm{im}(R_{i})$. We can show that these two spaces yield the same solution. 

\begin{lemma} 
\label{thm:samemin}
For $P_i$ and $R_i$ defined in \Cref{alg:orthomin(1)}, 
    \begin{align*}
         \argmin_{\rho \in \mathrm{im}(AP_{i-1}) \oplus \mathrm{im}(AR_i)} \| r^i - \rho\| = \argmin_{\rho \in \mathrm{im}(AP_{i})} \| r^i - \rho\|.
    \end{align*}
\end{lemma}
\begin{proof}
    Consider the problem 
    \begin{align*}
        \zeta^{i+1} = \argmin_{\zeta \in \Real^{N_I}} \| r^i - AP_i \zeta\|,
    \end{align*}
    yielding the $N_I$-ORTHOMIN(1) residual $r^{i+1} = r^i - AP_i\zeta^{i+1}$. If $\Pi_{AP_{i-1}}$ is the orthogonal projector on $\mathrm{im}(AP_{i-1})$, we can show that $\Pi_{AP_{i-1}} r^{i+1} = 0$. Indeed, 
    \begin{align*}
        \Pi_{AP_{i-1}}r^{i+1} = \Pi_{AP_{i-1}}r^i - \Pi_{AP_{i-1}}AP_i\zeta^{i+1}.
    \end{align*}
    However, $r^i \perp AP_{i-1}$ by Petrov-Galerkin orthogonality, so the first term is 0. Moreover, since $P_i$ is $A^\top A$-orthogonal to $P_{i-1}$, and 
    \begin{align*}
        \Pi_{AP_{i-1}} AP_i = AP_{i-1}(P_{i-1}^\top A^\top A P_{i-1})^{-1}\underbrace{P_{i-1}^\top A^\top A P_i}_{0},
    \end{align*}
    the second term is also 0. Consequently, the residual $r^{i+1}$ also satisfies the minimal residual problem with orthogonality constraint on $\mathrm{im}(AP_{i-1}) \oplus \mathrm{im}(AP_i)$, since it is orthogonal to $\mathrm{im}(AP_{i-1})$. Moreover, since the orthogonalization at line 8 gives $\mathrm{im}(R_{i}) \oplus \mathrm{im}(P_{i-1}) =\mathrm{im}(P_i) \oplus \mathrm{im}(P_{i-1})$, then $\mathrm{im}(AP_{i-1}) \oplus \mathrm{im}(AP_i) = \mathrm{im}(AP_{i-1}) \oplus \mathrm{im}(AR_i)$. By uniqueness of the minimizer, the statement of the lemma holds.
\end{proof}

This makes it clear that since $\mathrm{im}(R_i) = \mathcal{K}_{N_I}(A,r^i)$, by identifying $V_i = P_{i-1}$, we find that $N_I$-ORTHOMIN(1) is an augmented Krylov method that fits in the framework of \Cref{alg:augmentedpolynomialprecond}. Moreover, in GMRES-RB without any truncation due to SVD, we can show that the augmentation spaces $V_{i+1}$ are equivalently obtained by $A^\top A$-orthogonalization of $\mathcal{K}_{N_b}(A,r^i)$ against $V_i$, meaning the two algorithms are actually the same in exact arithmetic.

\begin{lemma}
\label{thm:sameaug}
    Assume that $AV_i$ has full column rank for any cycle $i$ of GMRES-RB($N_I$, $N_I$, $0$, $\infty$) and $N_I$-ORTHOMIN(1). Then, $\mathrm{im}(P_i) = \mathrm{im}(V_i)$ and $r^i_\text{GMRES-RB} = r^{i+1}_\text{ORTHOMIN}$ in exact arithmetic.
\end{lemma}
\begin{proof}
    To simplify the notation, the iterates, residuals and errors produced by GMRES-RB will have no subscript.
    
    After the warmup in GMRES-RB, the initial augmentation basis spans $\mathcal{K}_{N_I}(A,r^0_\mathrm{GMRES})$, the same as in $N_I$-ORTHOMIN(1), so $\mathrm{im}(V_0) = \mathrm{im}(P_0)$. Also at the end of the $N_I$ GMRES iterations in exact arithmetic, the final residual would be the same as for $N_I$-ORTHOMIN(1) since no restart happened yet, meaning $r^0=r^{N_b}_\text{GMRES}=r^1_\text{ORTHOMIN}$.

    During cycle $i\geq 0$, assume the induction hypotheses $\mathrm{im}(V_{i}) = \mathrm{im}(P_i)$ and $r^{i} = r^{i+1}_\text{ORTHOMIN}$. The GMRES-RB algorithm minimizes over the same space as $N_I$-ORTHOMIN(1):
    \begin{align*}
        \mathrm{im}(V_i) + \mathcal{K}_{N_I}(A,r^i) = \mathrm{im}(P_i) + \mathcal{K}_{N_I}(A,r_\text{ORTHOMIN}^{i+1}) = \mathrm{im}(P_i) \oplus \mathrm{im}(R_{i+1}).
    \end{align*}
    The first equality holds by the induction hypothesis, and \Cref{thm:samemin} tells us that the space on the right is the $N_I$-ORTHOMIN(1) space. Thus the residuals $r^{i+1}$ and $r^{i+2}_\text{ORTHOMIN}$ are equal.

The intermediate corrections problem, given by,
\begin{align*}
    \zeta^{k,i} = \argmin_{\zeta\in\Real^{k+m_i+1}} \|r^i - AW^{(k)}_i\zeta \|,
\end{align*}
can be rewritten as 
\begin{align*}  \zeta^{k,i}=\mqty[\alpha^{k,i}\\\beta^{k,i}],\quad \{\alpha^{k,i},\beta^{k,i}\}= \argmin_{\alpha,\beta} \|r^i - A[Z_i]_k\alpha - AV_i\beta \|.
\end{align*}
Thus, the intermediate correction can be written as $x^{k,i}-x^i = [Z_i]_k\alpha^{k,i} + V_i\beta^{k,i}$. For a fixed $\alpha$, the solution in $\beta$ is given by the normal equations
\begin{align*}
    \beta = (V_i^\top A^\top A V_i)^{-1}V_i^\top A^\top (r^i - A[Z_i]_{k}\alpha) = -(AV_i)^+ A[Z_i]_{k}\alpha,
\end{align*}
where $(\,\cdot\,)^+$ denotes the Moore-Penrose pseudo-inverse for matrices with full column rank. Since $\mathrm{im}(AV_i)$ is a subspace of $\mathrm{im}(AW_{i-1})$, and $r^i\perp \mathrm{im}(AW_{i-1})$, the first term is 0.
For the optimal $\alpha^{k,i}$, the correction is
\begin{align*}
    W_i\zeta^{k,i}_\circ = x^{k,i}-x^i = (I-V_i(AV_i)^+ A)[Z_i]_k\alpha^{k,i},
\end{align*}
which yields the following equalities, provided no breakdown happens:
\begin{align*}
    \mathrm{im}(V_{i+1})=\mathrm{im}(W_i[U_{i+1},\zeta^{i}])=\mathrm{im}(W_i[\zeta^{1,i}_\circ,\dots,\zeta^{i}])=(I-\Pi_{V_i}^{A^\top A})\mathcal{K}_{N_I}(A,r^i).
\end{align*}
Notice that $\Pi_{V_i}^{A^\top A}=V_i(AV_i)^+A$ is the $A^\top A$-orthogonal projector onto $\mathrm{im}(V_i)$, meaning that $V_{i+1}$ spans the Krylov subspace $\mathcal{K}_{N_I}(A,r^i)$, from which components of $\mathrm{im}(V_i)$ in the $A^\top A$-inner product have been removed. For $N_I$-ORTHOMIN(1), the updated subspace would be equal by using the induction hypothesis (notice that the projector is independent of the basis chosen)
\begin{align*}
    \mathrm{im}(P_{i+1})= (I-\Pi_{P_i}^{A^\top A}) \mathcal{K}_{N_I}(A,r_\text{ORTHOMIN}^{i+1}) = (I-\Pi_{V_i}^{A^\top A}) \mathcal{K}_{N_I}(A,r^{i}).
\end{align*}

By induction, since the initial subspaces and residuals are the same, they will remain the same since they are transformed in the same way.
\end{proof}

\begin{proposition}
\label{thm:equiv_orthomin}
The untruncated GMRES-RB($N_I$, $N_I$, $0$, $\infty$) method coincides, in exact arithmetic, with $N_I$-ORTHOMIN(1), i.e., $x^i_\textrm{GMRES-RB}=x^{i+1}_\mathrm{ORTHOMIN}$ for $i\geq 0$.
\end{proposition}
\begin{proof}
    Direct consequence of \Cref{thm:sameaug}. The shift by 1 comes from the fact that the iterate of the GMRES warmup phase is $x^0_\textrm{GMRES-RB}$, but is the first ORTHOMIN iterate.
\end{proof}

Now that this equivalence is established, we can use other equivalences proven by Chronopoulos \cite{Chronopoulos_1991}, for specific classes of matrices. In particular, when $A$ is symmetric, the method matches MINRES in exact arithmetic.
\begin{corollary}
    \label{thm:equiv_gmres}
    Let $M$ and $N$ be symmetric and skew-symmetric matrices respectively. Suppose $A=M$, $A=N$ or $A=I-N$. Moreover, assume that for some $N_I\geq 2$, $A^{N_I}$ is definite. Then, in exact arithmetic, GMRES and GMRES-RB($N_I$, $N_I$, $0$, $\infty$) produce the same iterates if they start from the same initial point $x^0$. More precisely, $x^i_\textrm{GMRES-RB} = x^{(i+1)N_I}_\textrm{GMRES}$ for $i\geq 0$.
\end{corollary}
\begin{proof}
    Theorem 5.5 and Corollary 5.1 of Chronopoulos \cite{Chronopoulos_1991} yield the result, because of an equivalence between $s$-ORTHOMIN(1) and $s$-GCR for some matrices, and then an equivalence between $s$-GCR and full GCR when $A^s$ is definite. GCR is equivalent to GMRES in exact arithmetic, see for instance the consequences of Lemma 6.21 from Saad \cite{Saad_2003}.
\end{proof}

This gives a first heuristic of why the proposed algorithm would work well. If for instance, the matrix $A$ is symmetric or skew-symmetric, and the truncated approximation space is a very good approximation of the untruncated one (e.g., there is fast decay of the singular values, by the Eckart-Young theorem \cite{Eckart_Young_1936}), then the method should approximate in some sense the behaviour of GMRES. It is not unlikely that the algorithm will also behave nicely when applied to matrices which are perturbations of symmetric or skew-symmetric matrices as a consequence. Moreover, in the case $N_I=2$ and where $A$ is symmetric, GMRES-RB(2,2,0,$\infty$) transfers two vectors in its reduced basis to the next cycle, meaning $W_i$ will contain four vectors, and that maintaining orthogonality to this window of four vectors is enough to construct the polynomial representation of the solution; GMRES-RB can capture the short-term recurrence if it exists.

For the GCRO-RB variant, we also see that in the case where all vectors are truncated, GCRO-RB($N_I+1$, $N_I$, 0, 0) is equivalent to GCROT($N_I$, 1), where only the final correction of the last cycle is kept. In a sense, GCROT($m$, $k$) is really the deflated analog of LGMRES($m$, $k$); they use the same information for the augmentation basis, but the former also projects the operator.

\subsection{Error Analysis}

For augmented Krylov methods, improved convergence rates or recycling properties are often obtained by using the hypothesis that the augmentation space approximates an appropriate $A$-invariant subspace \cite{Saad_1997,Eiermann_Ernst_Schneider_2000,Parks_de_Sturler_Mackey_Johnson_Maiti_2006,Kahl_Rittich_2017}. However, the hypothesis that the space generated by the snapshots contains or approximates well-enough the right invariant subspace cannot be verified a priori in general. Baker et al. \cite{Baker_Jessup_Manteuffel_2005,Baker_Jessup_Kolev_2009} propose convergence analyses based on angles between residuals to estimate how truncated GMRES or LGMRES oscillates. As mentioned in the previous section, we already capture the skip angles property because the last correction is in the augmentation subspace.

We study global behaviour of the GMRES-RB algorithm, meaning that we determine the error between GMRES-RB and full, untruncated GMRES.

\begin{lemma}
     The correction search spaces generated by GMRES-RB satisfy 
     \begin{align*}
         \mathrm{im}(V_{i}) &\subset \mathcal{K}_{N_b + iN_I}(A,r_\textrm{GMRES}^0),\\
          \mathrm{im}(W_{i}) &\subset \mathcal{K}_{N_b + (i+1)N_I}(A,r_\textrm{GMRES}^0),
     \end{align*}
     where $r^0_\textrm{GMRES}$ is the residual at the beginning of the warmup phase.
     \label{thm:inclusKrylov}
\end{lemma}
\begin{proof}

    After the warmup phase, $\mathrm{im}(V_0) \subset \mathcal{K}_{N_b}(A,r_\textrm{GMRES}^0)$. Also, $r^0 = r^{N_b}_\textrm{GMRES}\in\mathcal{K}_{N_b+1}(A,r_\textrm{GMRES}^0)$. Consequently, the first $N_I$ vectors of $W_0$, which span $\mathcal{K}_{N_I}(A,r^0)$, are all in $\mathcal{K}_{N_b + N_I}(A,r^0_\textrm{GMRES})$. The last $m_0$ columns are the columns of $V_0$, meaning that the whole $\mathrm{im}(W_0)$ is contained in $\mathcal{K}_{N_b + N_I}(A,r^0_\textrm{GMRES})$. The base case is satisfied; we assume the strong induction hypothesis
    \begin{align*}
         \mathrm{im}(V_{s}) &\subset \mathcal{K}_{N_b + sN_I}(A,r_\textrm{GMRES}^0),\\
          \mathrm{im}(W_{s}) &\subset \mathcal{K}_{N_b + (s+1)N_I}(A,r_\textrm{GMRES}^0),
    \end{align*}
    for $s\leq i$. By construction of the augmentation basis in \Cref{alg:gmresaccelere},
    \begin{align*}
        \mathrm{im}(V_{i+1}) \subset \mathrm{im}(W_i) \subset \mathcal{K}_{N_b + (i+1)N_I}(A,r_\textrm{GMRES}^0).
    \end{align*}
    Following the construction of $W_i$ in \Cref{alg:augmentedpolynomialprecond}, 
    \begin{align*}
        \mathrm{im}(W_{i+1}) = \mathrm{im}(V_{i+1}) + \mathcal{K}_{N_I}(A,r^{i+1}). 
    \end{align*}
    However, $r^{i+1}\in \mathcal{K}_{N_b+(i+1)N_I+1}(A,r^0_\mathrm{GMRES})$, since
    \begin{align*}
        r^{i+1} = r^0 - \sum_{j=0}^i AW_j\zeta^j,
    \end{align*}
    with $r^0\in \mathcal{K}_{N_b+1}(A,r_\textrm{GMRES}^0)$ and the second term being in $A\mathcal{K}_{N_b+(i+1)N_I}(A,r_\textrm{GMRES}^0)$. Thus, $\mathcal{K}_{N_I}(A,r^{i+1}) \subset \mathcal{K}_{N_b + (i+2)N_I}(A,r^0_\textrm{GMRES})$. To conclude,
    \begin{align*}
        \mathrm{im}(W_{i+1}) \subset \mathcal{K}_{N_b + (i+2)N_I}(A,r^0_\textrm{GMRES}).
    \end{align*}
\end{proof}  
This loose inclusion is sufficient to show the following proposition, which give the exact error terms of GMRES-RB when compared to full, untruncated GMRES. The accuracy is directly controlled by how well the vector $\Delta^i = x^{N_b+(i+1)N_I}_\textrm{GMRES} - x^i_\textrm{GMRES-RB}$ is approximated in the augmented subspace $\mathrm{im}(W_i)$ with respect to the $A^\top A$-inner product.

\begin{proposition}
\label{thm:PetrovGalerkinErr}
    For $i\geq 0$, let $e^{N_b+(i+1)N_I}_\textrm{GMRES}$ be the full GMRES error, and $e^{i+1}_\textrm{GMRES-RB}$ be the corresponding GMRES-RB error. They satisfy the following identity:
    \begin{align*}
        \| e^{i+1}_\textrm{GMRES-RB}\|_{A^\top A}^2 = \|e_\textrm{GMRES}^{N_b+(i+1)N_I}\|_{A^\top A}^2 + \|(I - \Pi^{A^\top A}_{W_i}) \Delta^i\|_{A^\top A}^2,
    \end{align*}
    where $\Pi^{A^\top A}_{W_i}$ is the orthogonal projector on $\mathrm{im}(W_i)$ w.r.t. the $A^\top A$-inner product.
\end{proposition}
\begin{proof}
    To shorten the notation in the proof, the GMRES-RB subscript will be dropped when referring to GMRES-RB iterates, residuals, and errors. Also, let $j:=N_b+(i+1)N_I$.
    
    Let $x^*$ be the exact solution to $Ax^*=b$. We know that $r^{i} = Ae^i$ where $e^i = x^* - x^i$. The residual satisfies
    \begin{align*}
        \| r^{i+1}\| = \min_{\zeta \in \Real^{N_I+m_i+1}} \| r^{i} - A W_i\zeta\|.
    \end{align*}
    Then, let $e^j_\textrm{GMRES} = x^* - x^j_\textrm{GMRES}$. We have that ${r_\textrm{GMRES}^j = Ae^j_\textrm{GMRES}}$ and by \Cref{thm:inclusKrylov}
    \begin{align*}
        e^i - e^j_\textrm{GMRES} = x^j_\textrm{GMRES} - x^i = (x^j_\textrm{GMRES}-x_\textrm{GMRES}^0) - \underbrace{(x^i-x_\textrm{GMRES}^0)}_{\mathcal{K}_j(A,r^0_\textrm{GMRES})} \in \mathcal{K}_j(A,r_\textrm{GMRES}^0).
    \end{align*}
    Then, we can additively shift the problem by $A(e^i - e^j_\textrm{GMRES})$,
    \begin{align*}
        \| r^{i+1}\|^2 = \min_{\zeta \in \Real^{N_I+m_i+1}} \|r_\textrm{GMRES}^j + A(e^i -e^j_\textrm{GMRES} - W_i\zeta)\|^2.
    \end{align*}
    The second term in the sum is in $A\mathcal{K}_{j}$ because of the inclusion in \Cref{thm:inclusKrylov}. Because $r_\textrm{GMRES}^j \perp A\mathcal{K}_{j}$, the Pythagorean theorem gives 
    \begin{align*}
        \| r^{i+1}\|^2 = \|r_\textrm{GMRES}^j\|^2 + \min_{\zeta\in\Real^{N_I+m_i+1}} \| A(x^j_\textrm{GMRES} - x^i - W_i\zeta)\|^2.
    \end{align*}
    By converting everything in terms of the error,
    \begin{align*}
        \| e^{i+1}\|_{A^\top A}^2 = \|e_\textrm{GMRES}^j\|_{A^\top A}^2 + \min_{w \in \mathrm{im}(W_i)} \| x^j_\textrm{GMRES} - x^i - w\|_{A^\top A}^2,
    \end{align*}
    and the last term is given by the orthogonal projector onto $\mathrm{im}(W_i)$ for the $A^\top A$-inner product.
\end{proof}

\Cref{thm:PetrovGalerkinErr} could be reinterpreted more generally in terms of geometry of nested Krylov subspaces. However, to keep it more understandable and interpretable, we will keep it specific to the case at hand and refer to the articles of Eiermann, Ernst and Schneider \cite{Eiermann_Ernst_Schneider_2000,Eiermann_Ernst_2001} for a more general study. Moreover, this conclusion is not unique to GMRES-RB, in the sense that the selection mechanism of the augmentation vectors did not matter in the proof. The only thing required is that the construction of \Cref{thm:inclusKrylov} remains true, i.e., that the augmentation vectors at the start of a cycle are subspaces of the right Krylov subspace. 

$x^i_\text{GMRES-RB}$ represents the starting point at beginning of cycle $i$. $\Delta^i$ represents the next full GMRES iterate, from which all current progress from GMRES-RB has been removed. Thus, the discrepancy between the next GMRES and GMRES-RB iterates is quantified by how much of the remaining necessary progress is orthogonal to the current search space, in the $A^\top A$-norm.

% Our analysis will aim to show that under the right conditions, the modes chosen by GMRES-RB can induce a per-cycle acceleration that would be absent from a simple restarted GMRES scheme.

% Let $\mathcal{W}_i = \mathrm{im}(V_i) \oplus \mathcal{K}_{N_I}(A,r^i)$.

% \begin{hypothesis}
%     Let $\widehat{r}^{i+1}$ be the residual obtained by GMRES($N_I$) starting from the residual $r^i$. Denote $\Pi_{A\mathcal{W}_i}$ the orthogonal projector onto $A\mathcal{W}_i$. We assume that, for some $\tau_i\in (0,1)$,
%     \begin{align*}
%         \frac{\| \Pi_{A\mathcal{W}_i} \widehat{r}^{i+1}\|}{\|\widehat{r}^{i+1}\|} > \tau_i.
%     \end{align*}
% \end{hypothesis}

% This means the closer $\tau_i$ close to 1, the better the information recycled will speed up the convergence within the cycle.

\section{Numerical Experiments}\label{sec:numerical}

% \subsection{Implementation}

% A simple framework will be used 

% Since \Cref{alg:gmresaccelere} shares the same structure as augmentation algorithms like LGMRES, and that GCRO-DR can be implemented with an arbitrary augmentation basis, the implementation of our GMRES-RB method will follow that of GCRO-DR, but by replacing harmonic Ritz vectors by POD vectors. 
% This introduces a significant difference in GCRO-DR and GMRES-RB, and it is that the intermediate GMRES corrections must be stored as snapshots, while only the final correction is computed at the end of a cycle in GCRO-DR.

% Incremental SVD \cite{Brand_2006} will be used to update directly the singular values and singular vectors everytime a new snapshot is added.

\subsection{Comparison to Other Solvers}

The reduced basis methods will be compared to classical,  state-of-the-art augmented Krylov methods. As a representative of methods of the GCRO family (\Cref{alg:augmentedgcromethod}), GCRO-DR, with either Rayleigh--Ritz or harmonic Ritz vectors will be used. This method will be implemented in the same framework as GCRO-RB, with only the augmented basis construction algorithm changing between the two. Methods which use corrections from previous cycles will also be compared against our proposed algorithms. The implementations of LGMRES and GCROT($m$,$k$) from SciPy 1.15.0 \cite{2020SciPy-NMeth} will be used. This implementation of GCROT($m$,$k$) also has the capability to use de Sturler's optimal truncation based on the largest principal angles, which will also be compared against. The baseline algorithms GMRES and GMRES($m$) will also be taken from SciPy. Also, the choice $N_b=N_I$ will be maintained to match the cycles of other methods.

\subsection{2D Convection-Diffusion Problem}

As a first example, we use the convection-diffusion problem \cite{Morgan_1995,Baker_Jessup_Manteuffel_2005}
\begin{align*}
    \Delta u + C\pdv{u}{x} = -(41)^2,
\end{align*}
for multiple $C$ values. The discretized operator is obtained from a centered finite difference stencil on the unit square, with $h=1/41$. The boundary conditions are homogeneous Dirichlet conditions. The initial value of the solver will be set to $x^0 = 0$.

For small $C$, the symmetric part of the discretized matrix is dominant with respect to the skew-symmetric part, and large $C$ reverses this relationship. \Cref{thm:equiv_gmres} gives the intuition that for symmetric or skew-symmetric matrices, the behaviour of the algorithm should become closer to that of GMRES when the truncation preserves most of the information. By sweeping $C$ values on this problem, we can observe the behaviour for matrices with a varying degree of dominance between the symmetric and skew-symmetric parts. The ratio of the norm of the skew-symmetric part to the norm of the symmetric part will vary between roughly 1/262 to 412, i.e., from mostly diffusive (symmetric) to mostly advective (skew-symmetric) problems. \Cref{tab:skewtosym} shows how the considered $C$ values and their associated ratio vary.

\begin{table}[H]
    \centering
    \caption{Ratio of the skew-symmetric part energy to symmetric part energy for multiple values of the advection coefficient $C$.}
    \begin{tabular}{c|c}
        $C$ &$\|A^\top - A\|_F/\|A^\top + A\|_F$\\
        \hline
        0 &  0\\
         1&  0.003817\\
         41& 0.1565\\
         1681& 6.417\\
         107901& 411.9\\
    \end{tabular}
    
    \label{tab:skewtosym}
\end{table}

\subsection{Comparison of the Reduced Basis Selection Rules}
\begin{figure}
    \centering
    \includegraphics[width=0.9\linewidth]{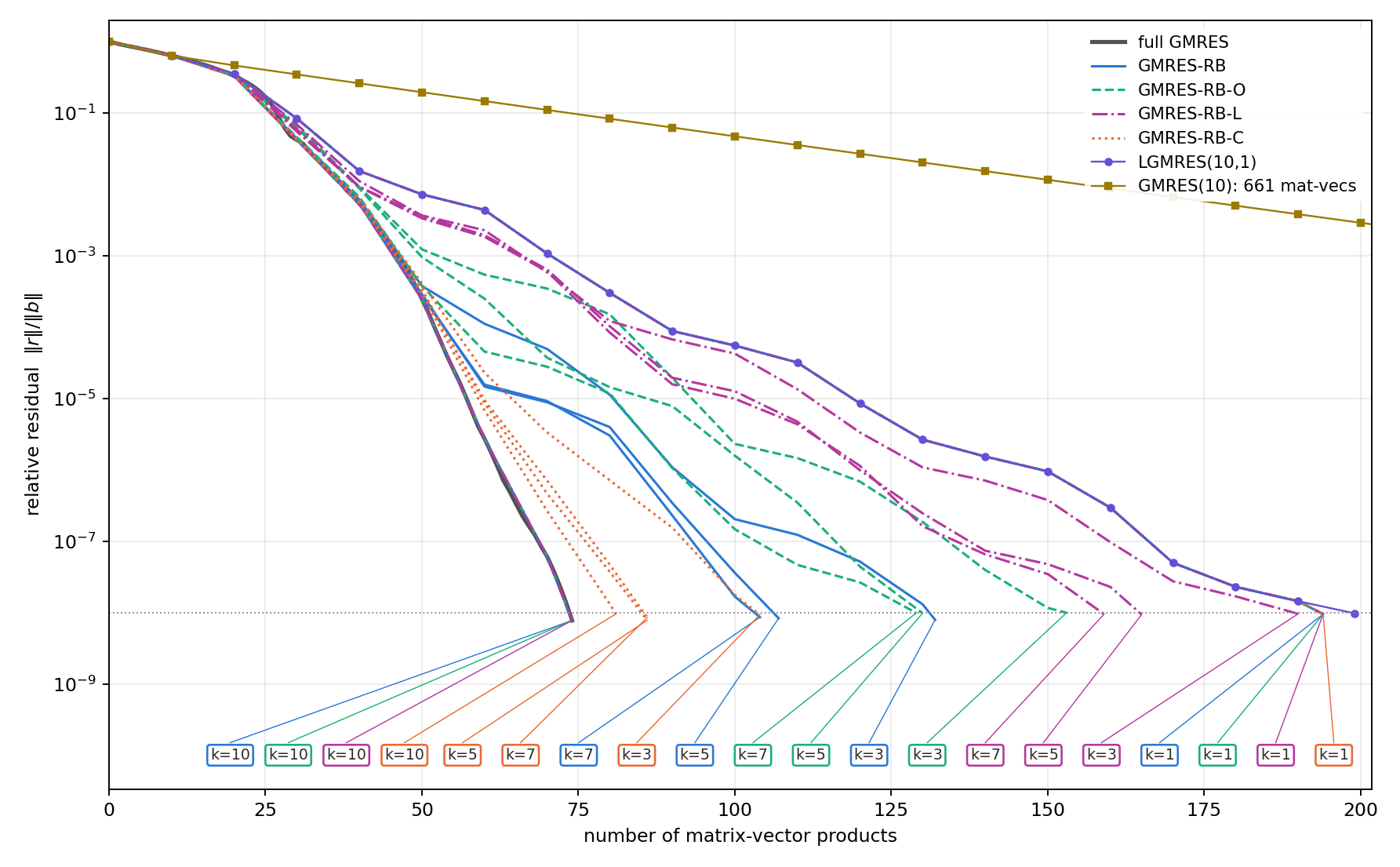}
    \caption{Comparison between the different variants of GMRES-RB at fixed budget $k$ on the Laplacian ($C=0$) problem, for $N_I=10$.}
    \label{fig:NI10COMPARE}
\end{figure}

\begin{figure}
    \centering
    \includegraphics[width=0.9\linewidth]{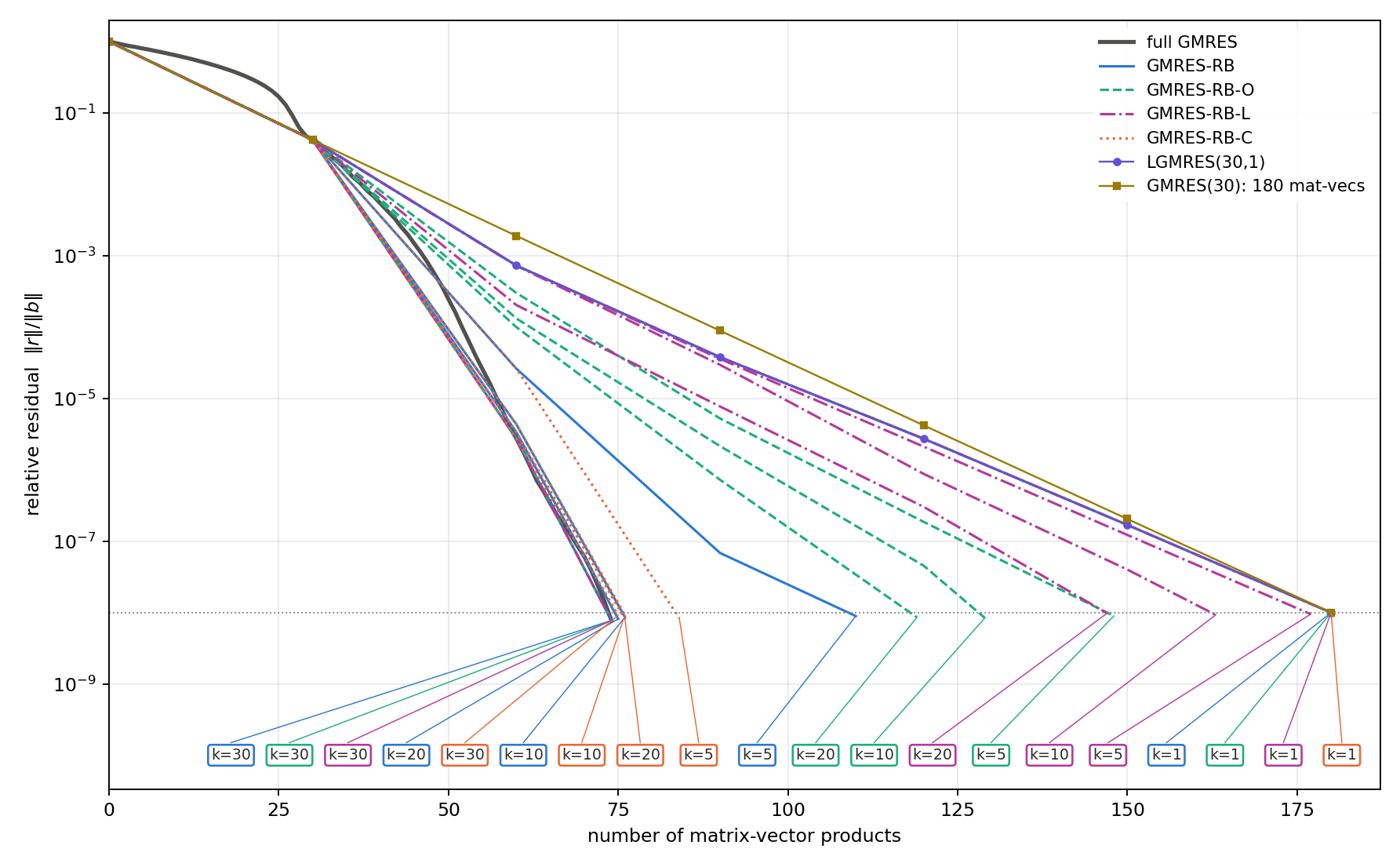}
    \caption{Comparison between the different variants of GMRES-RB at fixed budget $k$ on the Laplacian ($C=0$) problem, for $N_I=30$.}
    \label{fig:NI30COMPARE}
\end{figure}

The first test consists of comparing all the variants (-O, -L, -C, GCRO-RB) of GMRES-RB. To accurately measure which selection criteria/augmented subspace contains the richest information, all methods will be constrained to choose the same number of augmentation vectors, or the same budget, $k$. For the methods relying on an energy criterion (GMRES-RB and GMRES-RB-C), this will be enforced by setting $\varepsilon_\mathrm{rb}=0$ and $m_\mathrm{max}=k-1$, to keep the total correction and $k-1$ singular vectors. The basis size will be independent of $\varepsilon_\mathrm{rb}$, but the singular vectors will still be ranked according to their energy instead of their recency. In the case of GMRES-RB-C, the damping coefficient will be set to $\omega=0.5$. Moreover, we will validate that setting $k=N_I$ and $k=1$ indeed match the performance of full GMRES and LGMRES($N_I$, 1), by choosing $C=0$, a pure Laplace problem, which is SPD. \Cref{fig:NI10COMPARE,fig:NI30COMPARE} show the convergence of the relative residual norm for all variants at $k=1,3,5,7,10$ for $N_I=10$, and $k=1,5,10,20,30$ for $N_I=30$.

First, we notice that in both cases, full GMRES and LGMRES($N_I$,1) are exactly or almost exactly replicated by the budgets $k=N_I$ and $k=1$, due to round-off errors, which validates the conclusions of \Cref{thm:equiv_gmres}. Then, when comparing at equal budget, a ranking emerges. The best-performing methods at an equal budget are those based on the energy criterion; GMRES-RB-C first, then GMRES-RB. In the case $N_I=10$, we confirm that GMRES-RB-C is not equivalent to GMRES, even when there is no truncation, because the cumulative history prevents all corrections from the current cycle from being preserved. Thus, at $k=10$, it is outperformed by the other variants. However, as soon as truncation takes place, GMRES-RB-C requires fewer matvecs when compared with the others at equal budget. This means that information carried from multiple previous cycles does help to accelerate convergence. To measure how much information can be attributed to the historical augmented basis, we can measure the magnitude of the $j$-th column of $U_i$, $u^j$ in the recycled part of the basis, i.e., $\| R_{m_i}^{\top} u^j\|$. These values are averaged over multiple cycles in \Cref{tab:projection}. The most remarkable thing to notice is that the leading POD modes almost always lies exactly in the cumulative space, except for a few components. This is potentially due to the choice of $\omega=0.5$, which does not decrease sufficiently fast with respect to the magnitude of the snapshots through the cycles. However, as it will be seen in the following experiments, this choice gives a relatively performant method with respect to the other state-of-the-art methods that will be studied. The interpretation of this phenomenon is that the cumulative history compression process is a hybrid between picking the total correction of the last $k$ cycles, like LGMRES or GCROT would do, and the memoryless GMRES-RB, that only gets information from the current cycle. As will be seen in later experiments, although the exact composition of the subspace is not easily characterized, it proves to be, at times, richer than both the subspaces of LGMRES and GMRES-RB.

\begin{table}[htbp]
\centering
\caption{Component $\|R_{m_i}^\top u^j\|^2$ of the $k-1$ leading POD modes of GMRES-RB-C in the recycled space $V_i = W_{i-1}[U_i]_{m_i}$ (mean over modes, median over cycles), for $C = 0$, $N_I = 10$.}
\begin{tabular}{l|rrrrrrrrr}
 & $k=2$ & $k=3$ & $k=4$ & $k=5$ & $k=6$ & $k=7$ & $k=8$ & $k=9$ & $k=10$ \\
\hline
$\|R_{m_i}^\top u^j\|^2$ & 1.00 & 1.00 & 1.00 & 1.00 & 1.00 & 1.00 & 1.00 & 0.88 & 0.89 \\
\end{tabular}
\label{tab:projection}
\end{table}

In the $N_I=30$ case, we notice the high level of compressibility of the per-cycle history for this problem. Indeed, the SVD-based criterion in GMRES-RB and GMRES-RB-C, by keeping only $k=10$ of the 30 available snapshots (or $30+m_i$ for GMRES-RB-C), required only a handful more matvecs than GMRES needed to reach the same tolerance. In particular, for GMRES-RB-C, $k=5$ needed only around 15 more matvecs, and outperformed all the other variants by nearly 30 matvecs. 

It is also quite clear that for this problem, the information encoded by the oldest snapshots is richer than the latest, because at equal budget, GMRES-RB-O consistently requires 20 to 40 fewer matvecs than GMRES-RB-L. There are two possible hypotheses as to why this happens. First, the optimal case would be to keep all the snapshots. By measuring $\| (I-\Pi_{V_{i+1}}) S_i\|_F$, we can determine how much of the full space is captured by the reduced basis. \Cref{tab:projerror} shows that POD (which gives the optimal choice) is indeed lower than both other schemes, but the latest snapshots capture more of the full space at small $k$, which does not reflect the convergence behaviour observed in \Cref{fig:NI10COMPARE}, which is faster at any $k$.

\begin{table}[H]
\centering
\caption{Projection error $\| (I-\Pi_{V_{i+1}}) S_i\|_F/\|S_i\|_F$ of the full space on space spanned by the $k$ oldest/latest snapshots and leading POD modes, for $C = 0$ and $N_I = 10$.}
\begin{tabular}{l|ccccccccc}
Kept & $k=1$ & $k=2$ & $k=3$ & $k=4$ & $k=5$ & $k=6$ & $k=7$ & $k=8$ & $k=9$ \\
\hline
Oldest & 0.22 & 0.090 & 0.042 & 0.025 & 0.013 & 0.0072 & 0.0036 & 0.0018 & 0.00073 \\
Latest & 0.22 & 0.081 & 0.035 & 0.018 & 0.0091 & 0.0061 & 0.0040 & 0.0023 & 0.00094 \\
Energy & 0.22 & 0.053 & 0.0092 & 0.0045 & 0.0017 & 0.0011 & 0.00074 & 0.00050 & 0.00028 \\
\end{tabular}
\label{tab:projerror}
\end{table}

This leads to the second hypothesis, which is that the oldest snapshots span the Krylov space $\mathcal{K}_k(A,r^i)$, and the newest snapshots do not span a proper Krylov space, even if they are contained in one. In any case, the principal angles between the augmentation spaces spanned by the oldest/latest vectors and leading POD modes are shown in \Cref{tab:princangleO,tab:princangleL}. What it shows is that POD captures a large subspace of dimension $k-2$ or $k-3$ of both GMRES-RB-L and GMRES-RB-O. However, the extra dimensions are more aligned with the content of GMRES-RB-O. This hints at the fact that POD might approximate well a Krylov subspace with additional helpful directions.

\begin{table}[H]
\centering
\caption{Median principal angles $\theta_i$ (degrees) between the augmentation spaces built by GMRES-RB and GMRES-RB-O at budget $k$, for $C = 0$, $N_I = 10$.}
\begin{tabular}{c|rrrrrrrrr}
$k$ & $\theta_1$ & $\theta_2$ & $\theta_3$ & $\theta_4$ & $\theta_5$ & $\theta_6$ & $\theta_7$ & $\theta_8$ & $\theta_9$ \\
\hline
2 & 0 & 26.8 &  &  &  &  &  &  &  \\
3 & 0 & 1.8 & 48.0 &  &  &  &  &  &  \\
4 & 0 & 0.1 & 10.0 & 69.6 &  &  &  &  &  \\
5 & 0 & 0 & 0.7 & 19.9 & 79.1 &  &  &  &  \\
6 & 0 & 0 & 0 & 1.0 & 26.3 & 81.7 &  &  &  \\
7 & 0 & 0 & 0 & 0 & 0.8 & 20.8 & 79.5 &  &  \\
8 & 0 & 0 & 0 & 0 & 0 & 0 & 10.9 & 68.7 &  \\
9 & 0 & 0 & 0 & 0 & 0 & 0 & 0 & 0 & 48.0 \\
\end{tabular}
\label{tab:princangleO}
\end{table}

\begin{table}[htbp]
\centering
\caption{Median principal angles $\theta_i$ (degrees) between the augmentation spaces built by GMRES-RB and GMRES-RB-L at budget $k$, for $C = 0$, $N_I = 10$.}
\begin{tabular}{c|rrrrrrrrr}
$k$ & $\theta_1$ & $\theta_2$ & $\theta_3$ & $\theta_4$ & $\theta_5$ & $\theta_6$ & $\theta_7$ & $\theta_8$ & $\theta_9$ \\
\hline
2 & 0 & 31.7 &  &  &  &  &  &  &  \\
3 & 0 & 5.3 & 62.9 &  &  &  &  &  &  \\
4 & 0 & 0.7 & 26.5 & 84.3 &  &  &  &  &  \\
5 & 0 & 0 & 4.0 & 50.1 & 86.8 &  &  &  &  \\
6 & 0 & 0 & 0.1 & 9.1 & 58.5 & 88.3 &  &  &  \\
7 & 0 & 0 & 0 & 0 & 8.4 & 61.1 & 88.0 &  &  \\
8 & 0 & 0 & 0 & 0 & 0 & 0 & 47.9 & 86.8 &  \\
9 & 0 & 0 & 0 & 0 & 0 & 0 & 0 & 0 & 82.2 \\
\end{tabular}
\label{tab:princangleL}
\end{table}

\subsection{The Diffusion-Dominated Case}
\begin{figure}
    \centering
    \includegraphics[width=0.9\linewidth]{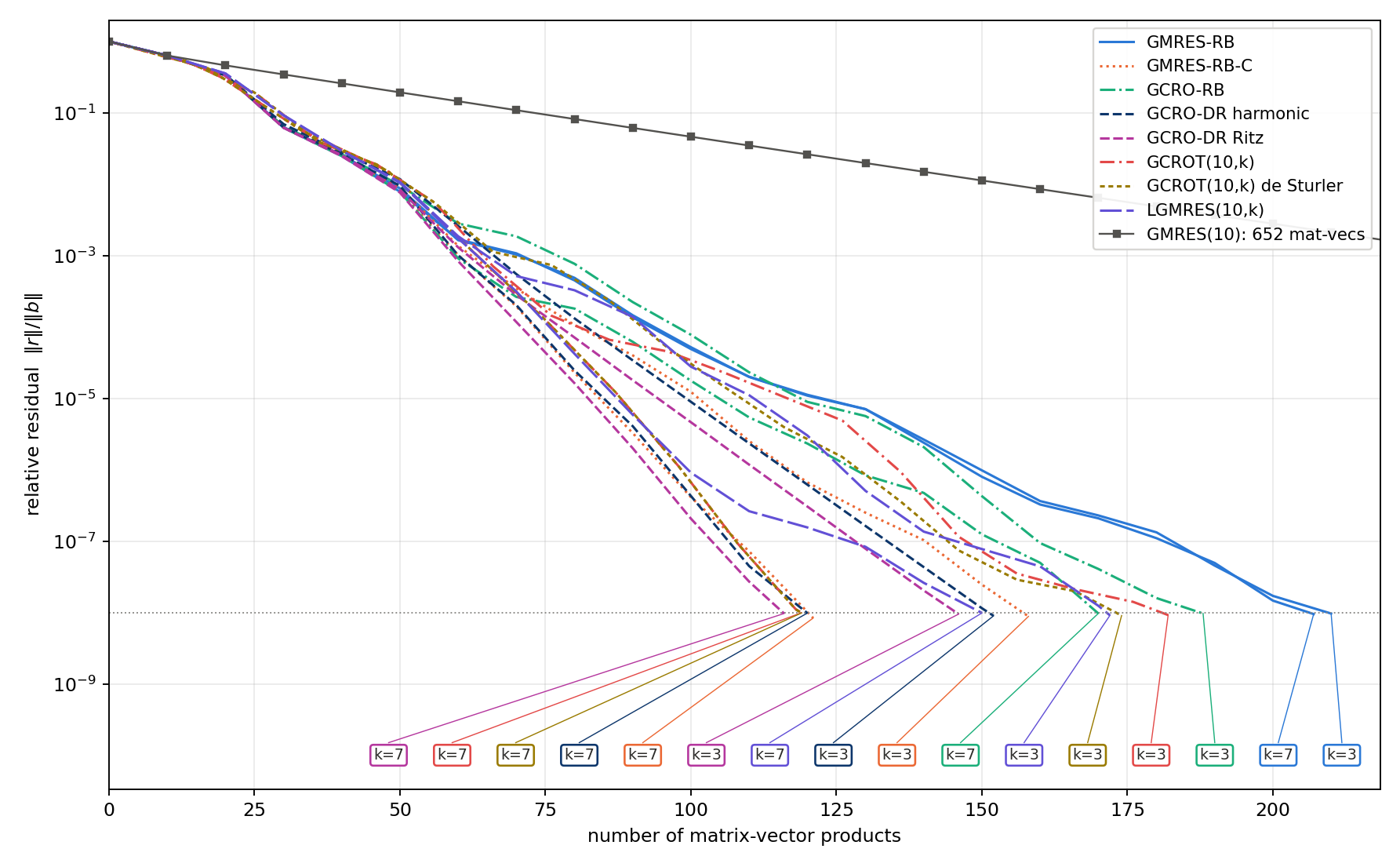}
    \caption{Comparison between variants of GMRES-RB and other classical solvers on the convection-diffusion problem ($C=1$) with a fixed budget $k$, and $N_I=10$.}
    \label{fig:NI10_cd_1}
\end{figure}

\begin{figure}
    \centering
    \includegraphics[width=0.9\linewidth]{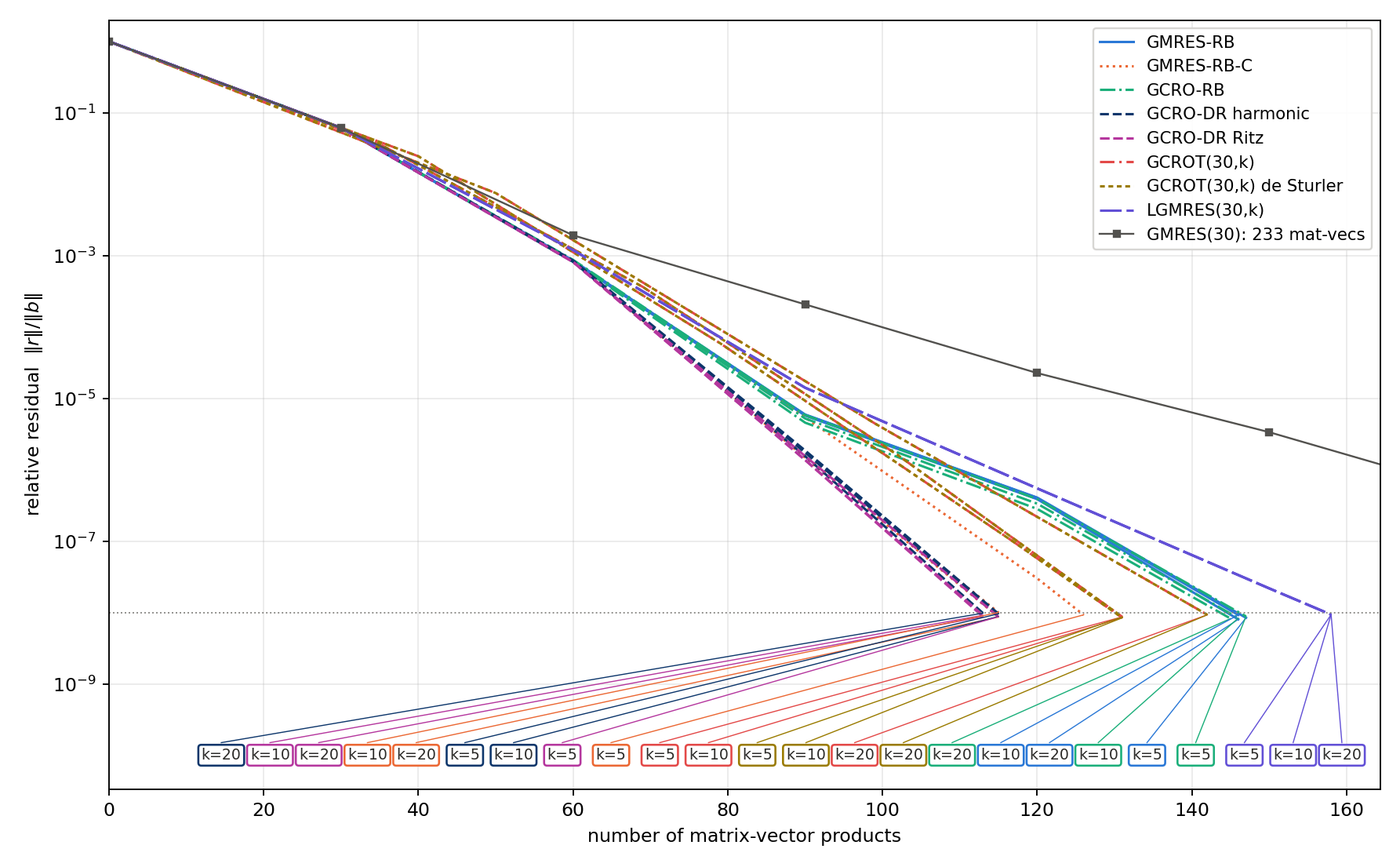}
    \caption{Comparison between variants of GMRES-RB and other classical solvers on the convection-diffusion problem ($C=1$) with a fixed budget $k$, and $N_I=30$.}
    \label{fig:NI30_cd_1}
\end{figure}

\begin{figure}
    \centering
    \includegraphics[width=0.9\linewidth]{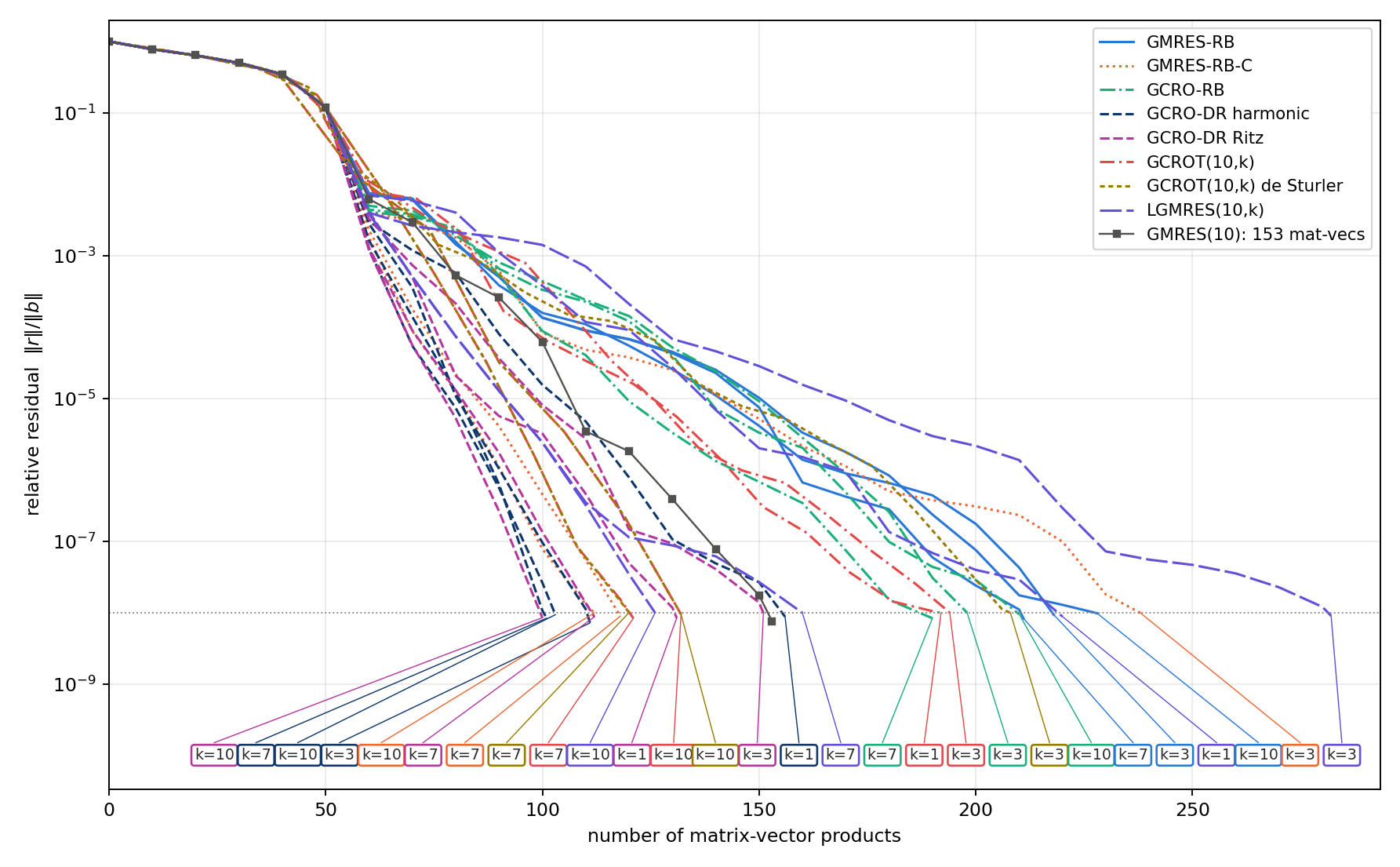}
    \caption{Comparison between variants of GMRES-RB and other classical solvers on the convection-diffusion problem ($C=41$) with a fixed budget $k$, and $N_I=10$.}
    \label{fig:NI10_cd_41}
\end{figure}

\begin{figure}
    \centering
    \includegraphics[width=0.9\linewidth]{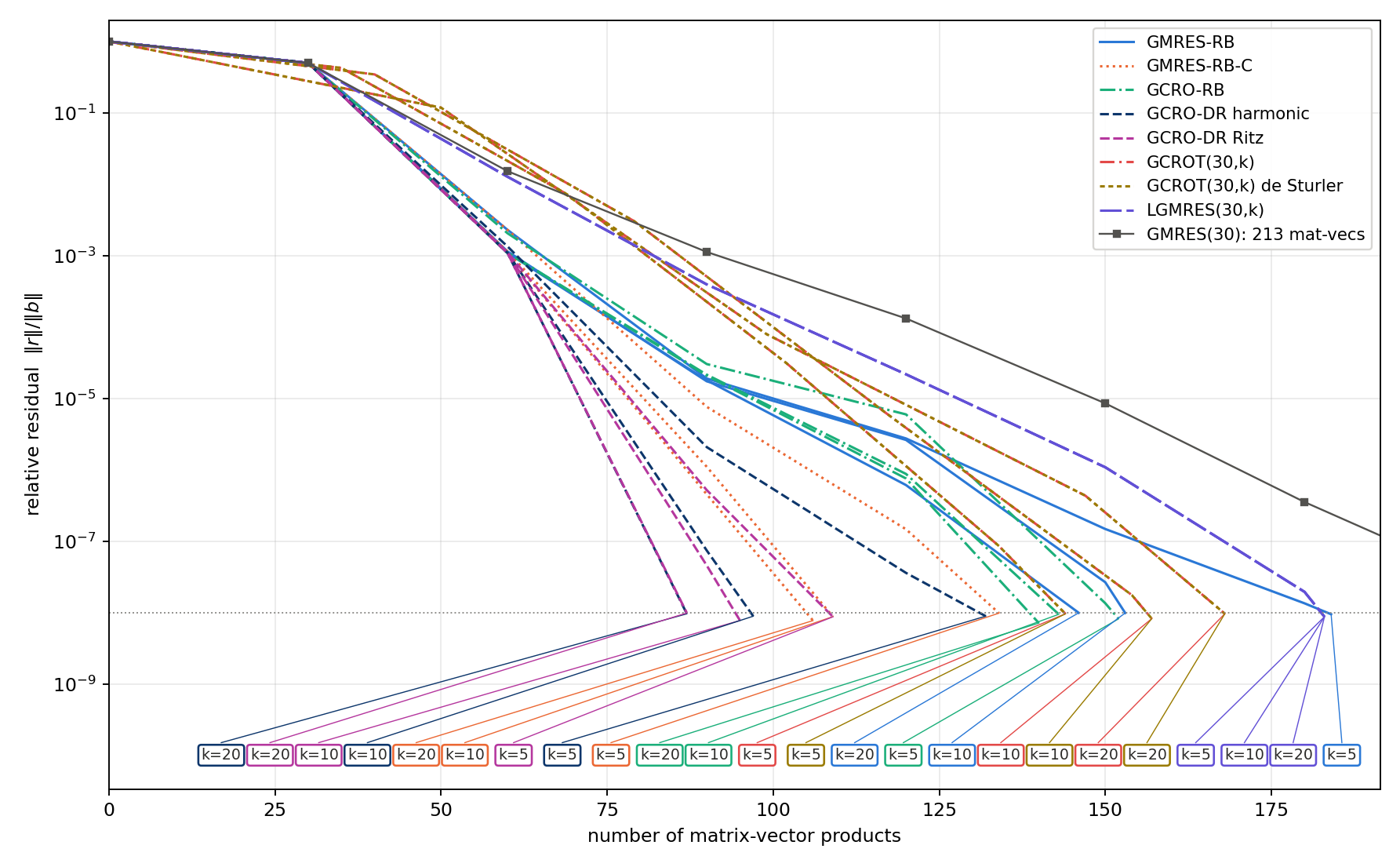}
    \caption{Comparison between variants of GMRES-RB and other classical solvers on the convection-diffusion problem ($C=41$) with a fixed budget $k$, and $N_I=30$.}
    \label{fig:NI30_cd_41}
\end{figure}

The cases where $C=1$ and $C=41$ will be studied. For these values, the operator is still diffusion-dominated. \Cref{fig:NI10_cd_1} shows convergence curves for $N_I=10$ and \Cref{fig:NI30_cd_1} for $N_I=30$ when $C=1$ (resp. \Cref{fig:NI10_cd_41,fig:NI30_cd_41} in the case $C=41$). We note that since $k=1$ is equivalent to LGMRES for GMRES-RB and GMRES-RB-C, they will be excluded from the figures if $k=1$ is present.

The most striking feature of these figures is that at equal budgets, the deflation methods GCRO-DR with either harmonic Ritz or Ritz vectors outperform all other methods. This is to be expected, since these methods approximate invariant subspaces, and benefit from provably improved convergence rates.

One of the noticeable things that can be observed is that in the case $N_I=10$, $C=41$, the behaviour in $k$ is nonmonotonic. For multiple methods, notably GCRO-DR with harmonic vectors, GMRES-RB and GCRO-RB, $k=7$ performs significantly better than $k=10$.

The plain GMRES-RB method is generally lagging behind every other algorithm at equal budget, except sometimes for LGMRES at higher $C$ and $N_I$. This behaviour is slightly surprising, since in the $C=0$ case, GMRES-RB largely outperformed GCRO-RB at the same budget, but introducing even a relatively small amount of advection ($C=1$, skew-sym. to sym. ratio well under 1\%) reversed the trend. This indicates that the behaviour of GMRES-RB is very sensitive to departure from symmetry. On the contrary, GMRES-RB-C, for largest $k$ values, seems to outperform or at least approximately match both GCROT variants at equal budget. However, there seems to be some sensitivity to the value of $k$, since in the case $C=41$ and $N_I=10$, taking $k=3$ leads GMRES-RB-C to be the second worst method. What is also surprising is that GMRES(10) outperforms most reduced basis methods and classical methods for $C=41$ and $N_I=10$.

\newpage
\subsection{The Advection-Dominated Case}

\begin{figure}[htbp]
    \centering
    \includegraphics[width=0.9\linewidth]{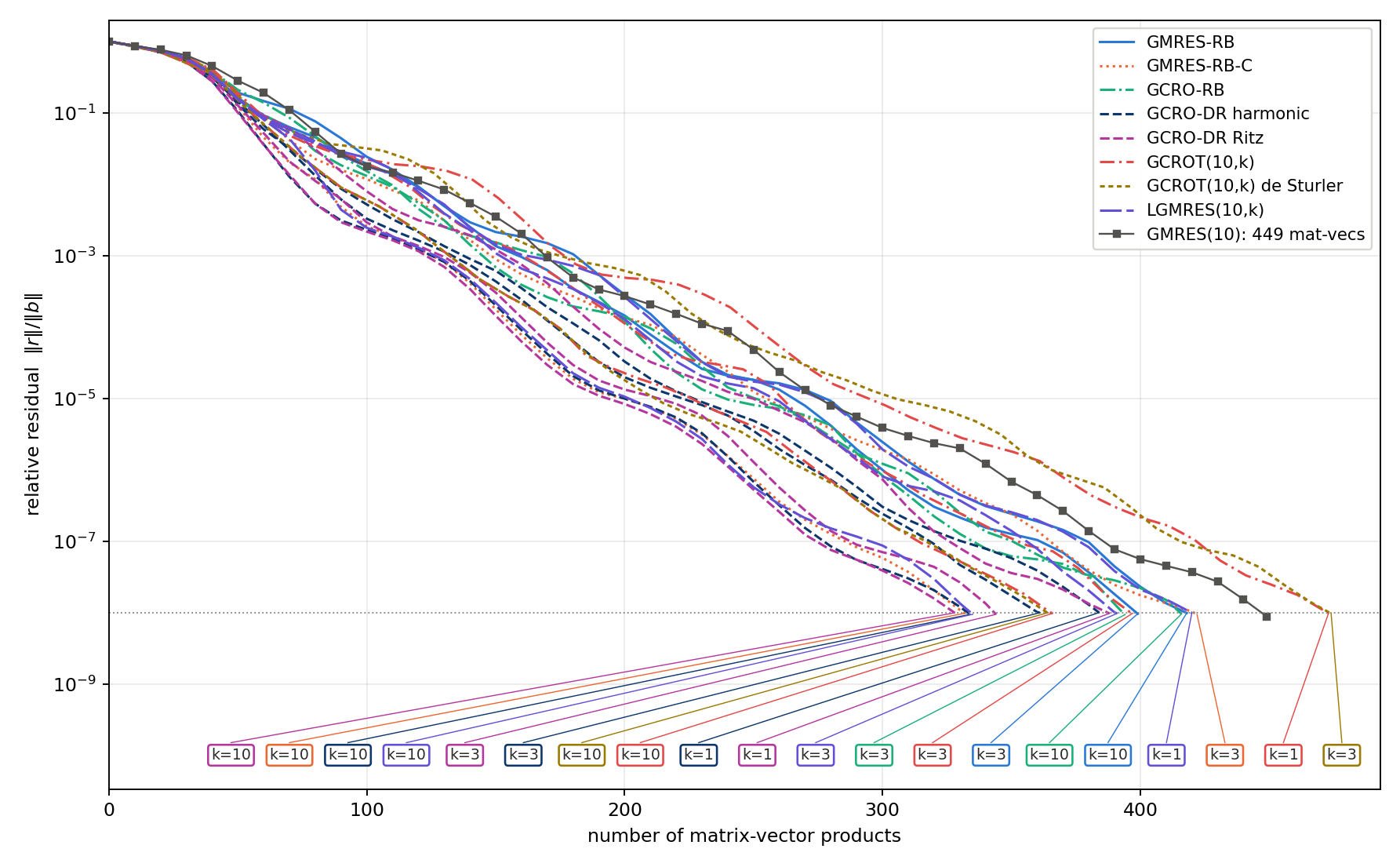}
    \caption{Comparison between variants of GMRES-RB and other classical solvers on the convection-diffusion problem ($C=1681$) with a fixed budget $k$, and $N_I=10$.}
    \label{fig:NI10_cd_1681}
\end{figure}

\begin{figure}[htbp]
    \centering
    \includegraphics[width=0.9\linewidth]{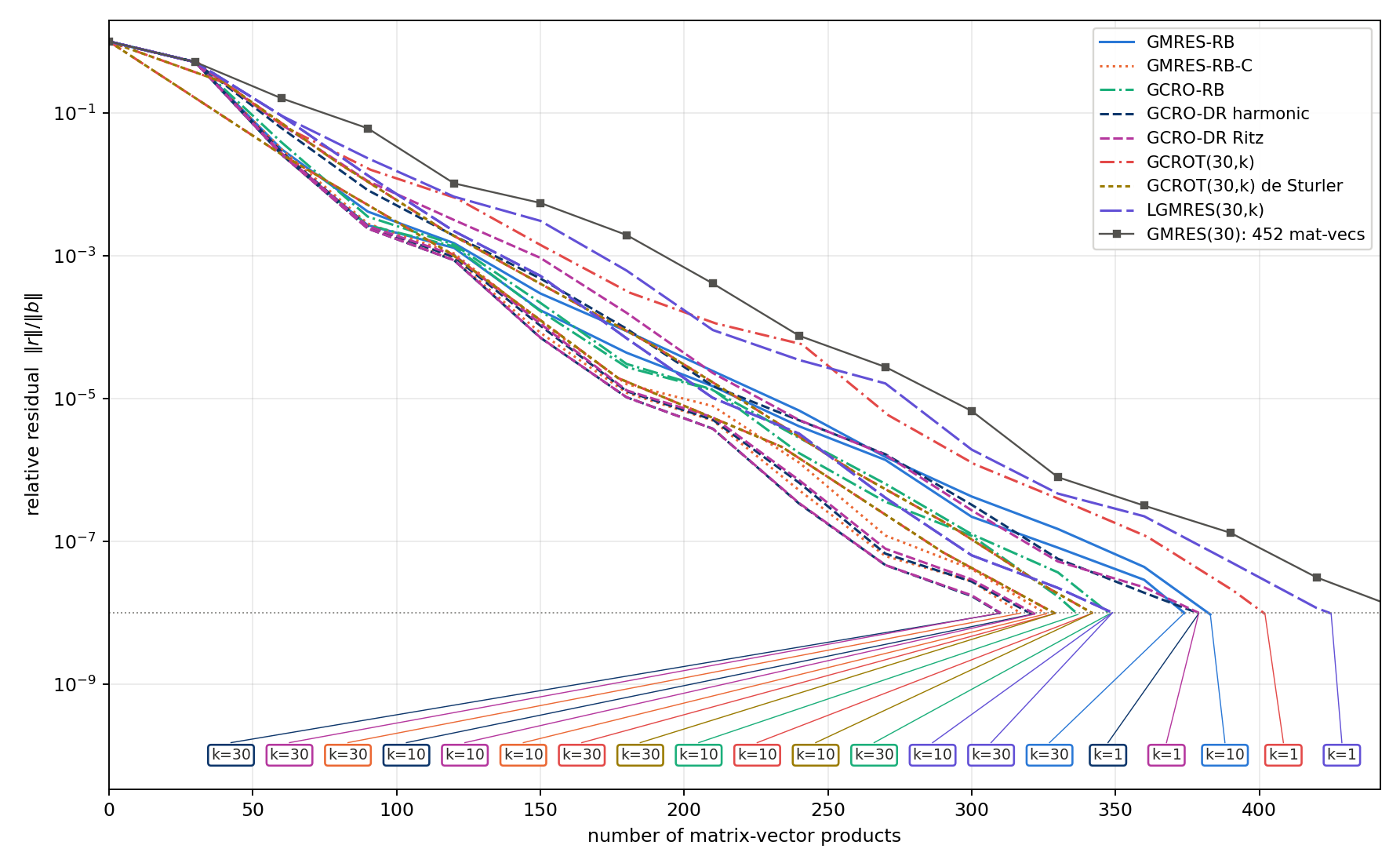}
    \caption{Comparison between variants of GMRES-RB and other classical solvers on the convection-diffusion problem ($C=1681$) with a fixed budget $k$, and $N_I=30$.}
    \label{fig:NI30_cd_1681}
\end{figure}

The case $C=1681$ yields a ratio of approximately 6 between the skew-symmetric and symmetric parts. \Cref{fig:NI10_cd_1681} shows what happens at $N_I=10$, and \Cref{fig:NI30_cd_1681} at $N_I=30$. In that case, GMRES-RB-C again outperforms the other reduced-basis methods. In fact, in this regime, the method performs almost as well as (or even slightly outperforms, in the case $N_I=10$) the deflation methods at larger $k$. Most of the time, it beats the tested GCROT methods. However, at small $k$, the trend reverses and it is almost always the worst method. As for GCRO-RB and GMRES-RB the former still requires fewer matvecs than the latter, especially at $N_I=30$.

\begin{figure}[htbp]
    \centering
    \includegraphics[width=0.9\linewidth]{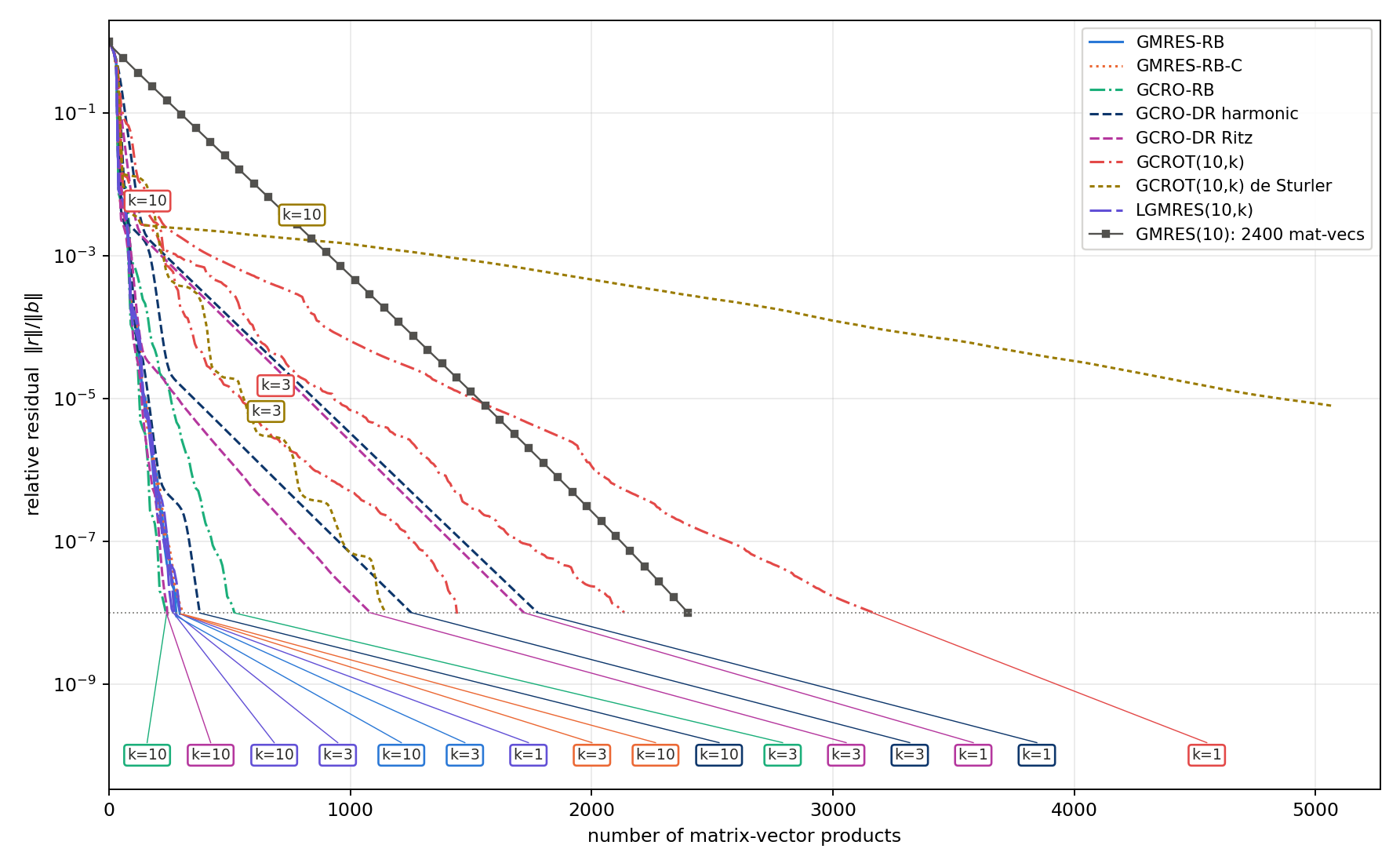}
    \caption{Comparison between variants of GMRES-RB and other classical solvers on the convection-diffusion problem ($C=107901$) with a fixed budget $k$, and $N_I=10$.}
    \label{fig:NI10_cd_107901}
\end{figure}

\begin{figure}[htbp]
    \centering
    \includegraphics[width=0.9\linewidth]{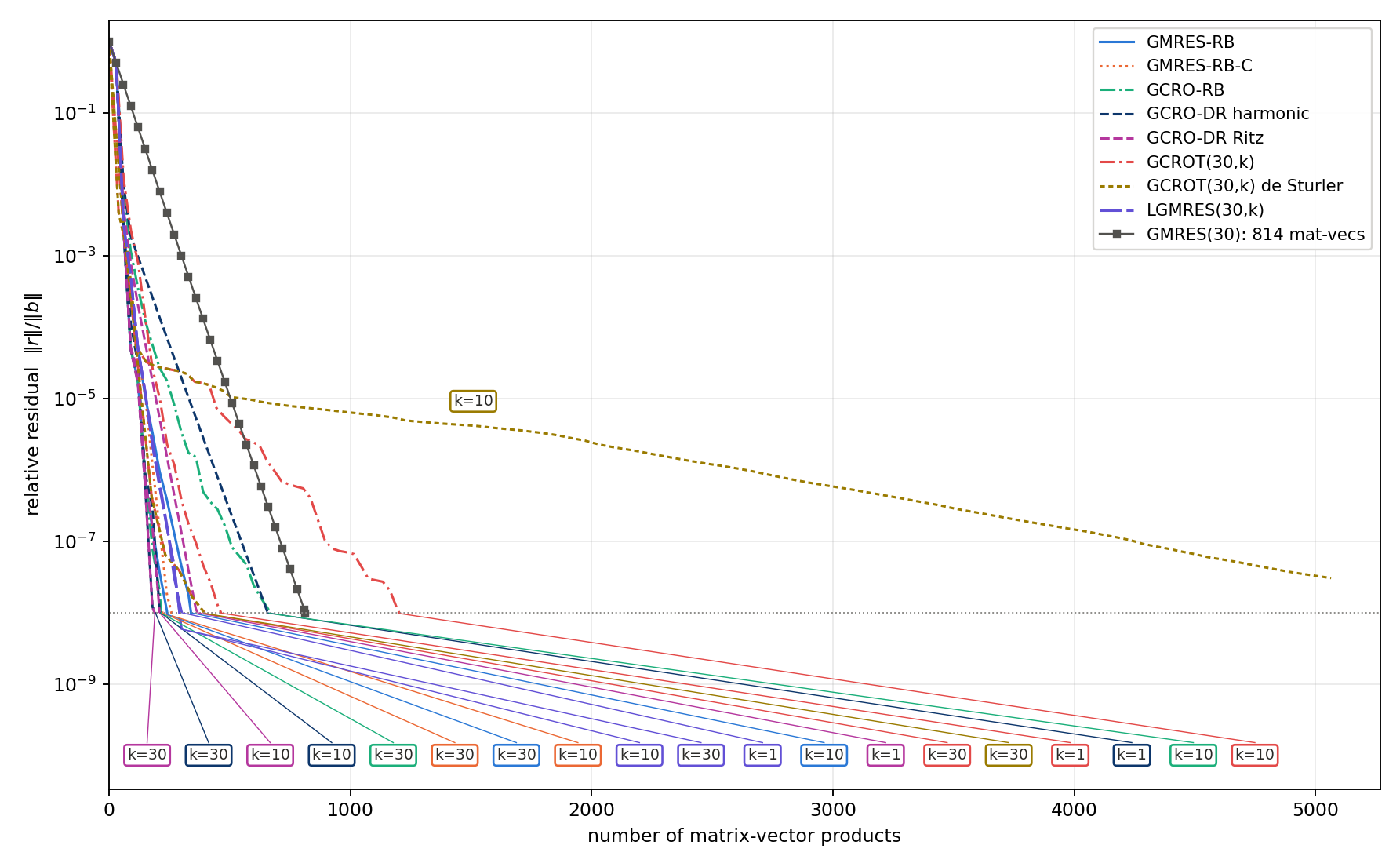}
    \caption{Comparison between variants of GMRES-RB and other classical solvers on the convection-diffusion problem ($C=107901$) with a fixed budget $k$, and $N_I=30$.}
    \label{fig:NI30_cd_107901}
\end{figure}

The last parameter value studied, $C=107901$, makes the problem largely dominated by advection, and the relative distance from $A$ to its skew-symmetric part is around $0.2\%$. \Cref{fig:NI10_cd_107901,fig:NI30_cd_107901} show the convergence behaviour at $N_I=10$ and $N_I=30$ respectively. Because well-performing methods are clustered between 100 and 400 matvecs, the number of matvecs for each budget is also shown in \Cref{tab:NI10_largeC,tab:NI30_largeC}.

\begin{table}[htbp]
\centering
\caption{Number of matvecs required to reach convergence, $C = 107901$, $N_I = 10$. The number in bold is the lowest in its column.}
\begin{tabular}{l|cccccc}
Method & $k=1$ & $k=2$ & $k=3$ & $k=5$ & $k=7$ & $k=10$ \\
\hline
GMRES-RB & \textbf{294} & \textbf{270} & 294 & \textbf{290} & 290 & 280 \\
GMRES-RB-C & \textbf{294} & 284 & 298 & 348 & 286 & 308 \\
GCRO-RB & 2692 & 318 & 520 & 404 & 332 & \textbf{238} \\
GCRO-DR harmonic & 1776 & 1776 & 1252 & 746 & 418 & 376 \\
GCRO-DR Ritz & 1720 & 1720 & 1084 & 548 & \textbf{252} & 242 \\
GCROT & 3163 & 439 & 1442 & 1241 & 1912 & 2135 \\
GCROT de Sturler & 3163 & 429 & 1148 & -- & -- & -- \\
LGMRES & 296 & 322 & \textbf{274} & \textbf{290} & 288 & 266 \\
\end{tabular}
\label{tab:NI10_largeC}
\end{table}

\begin{table}[htbp]
\centering
\caption{Number of matvecs required to reach convergence, $C = 107901$, $N_I = 30$. The number in bold is the lowest in its column.}
\begin{tabular}{l|cccccccc}
Method & $k=1$ & $k=2$ & $k=3$ & $k=5$ & $k=7$ & $k=10$ & $k=20$ & $k=30$ \\
\hline
GMRES-RB & \textbf{302} & \textbf{274} & 264 & 300 & 312 & 340 & 350 & 242 \\
GMRES-RB-C & \textbf{302} & 296 & 278 & \textbf{238} & 272 & 262 & 268 & 222 \\
GCRO-RB & 432 & 494 & 558 & 752 & 702 & 664 & 402 & 216 \\
GCRO-DR harmonic & 658 & 658 & \textbf{260} & 268 & 216 & 214 & \textbf{198} & 192 \\
GCRO-DR Ritz & 368 & 368 & 280 & 244 & \textbf{214} & \textbf{210} & \textbf{198} & \textbf{190} \\
GCROT & 465 & 965 & 1030 & 1212 & 1594 & 1205 & 297 & 398 \\
GCROT de Sturler & 465 & 941 & 870 & -- & -- & -- & 297 & 398 \\
LGMRES & \textbf{302} & 326 & 312 & 300 & 300 & 300 & 300 & 300 \\
\end{tabular}
\label{tab:NI30_largeC}
\end{table}

The global observation here is that for this problem setup, GMRES-RB starts to show some promise again. On larger budgets $k\geq 7$, GCRO-DR with Rayleigh--Ritz vectors starts picking up and it significantly reduces the number of matvecs required. For $N_I=k=10$, GCRO-RB is also able to beat the deflation method by a few matvecs. However, for smaller budgets $k\leq 5$, GMRES-RB and GMRES-RB-C always come on top, except once to GCRO-DR and once to LGMRES (as mentioned earlier, in the case $k=1$, it is equivalent to LGMRES in exact arithmetic). The magnitude of the gains when compared to classical methods is modest in most cases, but the best recorded case happens at budget $k=2$. The closest method to GMRES-RB is LGMRES for both $N_I$ values, and it requires 19\% more matvecs than the former.

\section{Conclusion}

This work introduces a family of accelerated Krylov methods, based on the framework of on-the-fly reduced bases. We show that these methods fit naturally in the framework of augmented Krylov methods, the main difference being that the augmentation vectors are obtained from a data-driven criterion.

We show that if the basis is not truncated, i.e., it consists of all intermediate corrections completed in the current cycle, the proposed GMRES-RB algorithm is equivalent in arithmetic to an $s$-step ORTHOMIN($k$) method. By using known results about this algorithm, it is possible to deduce that for symmetric or skew-symmetric matrices, the existence of a short-recurrence orthogonalization algorithm will ensure that GMRES-RB also coincides with GMRES in exact arithmetic. We can conjecture that matrices that are either close to symmetric or skew-symmetric could make the reduced basis methods behave well.

We then illustrate numerically that some methods in the proposed family, especially GMRES-RB-C, can perform almost as well as deflation methods for a similar budget. It also often outperforms correction-based augmentation methods like GCROT and LGMRES, especially when the matrices are very close to being either symmetric or skew-symmetric, and for larger $k$ values. The plain GMRES-RB method also performs well in these extremal cases, but for general, nonsymmetric matrices, it is generally outperformed by other augmentation algorithms. This allows us to conclude that preserving some convergence history from previous cycles as well as from the current cycle is important to preserve an improved global convergence rate, except for instance when the matrix is SPD, where GMRES is matched exactly by untruncated GMRES-RB.

A follow-up problem to the conclusions of this work is to establish exactly the rate at which departure from symmetry or skew-symmetry degrades the performance of GMRES-RB. Another relevant question, in the case of GMRES-RB-C, is to understand how departure from symmetry is compensated by the presence of the compressed convergence history, and how to balance in a more refined way the choice between information from older and later cycles. In particular, we hope that this could bring greater stability to GMRES-RB-C, because smaller $k$ values made it behave poorly.  Moreover, it would be interesting to relate the properties of $A$ and the decay of the singular values of the snapshot matrix, which would tell us how many modes can be truncated in each situation without degrading too much the short-recurrence in the symmetric case. More generally, the understanding of the mechanisms behind the improved convergence with respect to GMRES($m$), even in the SPD case, needs to be improved.

\section*{Acknowledgements}
The author thanks Prof. F\'elix Kwok and Prof. Jean Deteix for enlightening discussions and advice about the ideas and content of this manuscript.

\section*{Conflict of Interest}

The author declares no conflict of interest.

\section*{Funding}

The author was funded by a Natural Sciences and Engineering Research Council of Canada (NSERC) grant (grant number ALLRP 585256-2023).

\section*{AI Usage Statement}

The author used the large language model (LLM) Claude Opus 5.5 in the research leading up to the writing of this article, between September and November 2026. It was used for proof-reading the text and theoretical results of this manuscript, for code generation, optimization and debugging, and more generally as a conversational agent for testing intuitions and exchanging proof ideas. No text or data directly output by the LLM was included in the manuscript. Every conclusion presented in this work that was inspired by an output of the LLM was thoroughly validated and edited by the author.

\bibliography{ref}

\end{document}